# N-Adaptive Ritz Method: A Neural Network Enriched Partition of Unity for Boundary Value Problems

## Jonghyuk Baek[1,2], Yanran Wang[1], and J. S. Chen[1*]

[1] Department of Structural Engineering, University of California, San Diego, La Jolla, CA 92093, USA

[2] Coreform, LLC, Orem, UT 84097, USA

## Abstract

Conventional finite element methods are known to be tedious in adaptive refinements due to their conformal regularity requirements. Further, the enrichment functions for adaptive refinements are often not readily available in general applications. This work introduces a novel neural network-enriched Partition of Unity (NN-PU) approach for solving boundary value problems via artificial neural networks with a potential energy-based loss function minimization. The flexibility and adaptivity of the NN function space are utilized to capture complex solution patterns that the conventional Galerkin methods fail to capture. The NN enrichment is constructed by combining pre-trained feature-encoded NN blocks with an additional untrained NN block. The pre-trained NN blocks learn specific local features during the offline stage, enabling efficient enrichment of the approximation space during the online stage through the Ritz-type energy minimization. The NN enrichment is introduced under the Partition of Unity (PU) framework, ensuring convergence of the proposed method. The proposed NN-PU approximation and feature-encoded transfer learning forms an adaptive approximation framework, termed the neural-refinement (n-refinement), for solving boundary value problems. Demonstrated by solving various elasticity problems, the proposed method offers accurate solutions while notably reducing the computational cost compared to the conventional adaptive refinement in the mesh-based methods.

 * Corresponding author.

 *E-mail address*: jsc137@ucsd.edu (J.S. Chen).

# 1 Introduction

The adaptability and robust representation capabilities of neural networks (NNs) [1, 2] have propelled advancements in computational mechanics, contributing notably to various aspects such as data-driven constitutive modeling [3-8], data-driven computational mechanics [9-15], multiscale modeling [16-20], and damage and fracture modeling [21, 22]. Leveraging their inherent adaptivity, NNs serve as effective function approximators for solving general partial differential equations (PDEs). Their ability to construct flexible function spaces makes NNs an appealing alternative to the mesh-based method, such as the generalized finite element method [23] and the global-local enrichment [24], particularly for problems that involve localized features.

As data accessibility grows and platforms like Theano [25], TensorFlow [26], MXNet [27], and Keras [28] offer features like high-performance computing and automatic differentiation [29], neural networks emerge as a new approach to solve mathematical models. Raissi et al. [30] introduced physics-informed neural networks (PINNs), using prior PDE knowledge to estimate physical systems with sparse data via residual-based minimization based on collocation framework. Drawing from the concept of PINNs, Haghighat and Juanes [31] created SciANN, a Python wrapper tailored for scientific computing, successfully capturing localization behaviors in perfectly plastic materials. Sirignano and Spiliopoulos [32] combined Galerkin methods with deep learning, leveraging deep neural networks for function approximation and minimizing residual-based objective functions at random collocation points. However, a major challenge of using NNs for PDE solved via residual-based loss functions with strong form collocation is the extensive computational cost, requiring large amount of sampling points across the problem domain for adequate solution resolution. An alternative method is formulated with variational problems, where the corresponding functional acts as an energy descriptor for mechanical systems, serving as a loss function to optimize NNs. Weinan and Yu [33] introduced the Deep Ritz method, employing deep neural networks as trial functions for solving PDEs by minimizing energy, adaptable to higher dimensions. However, the resulting variational problems can be non-convex, leading to challenges like local minima and saddle points, and imposing essential boundary conditions presents non-trivial complexities. Samaniego et al. [34] introduced NNs as function approximation machines and use the energy of mechanical system as a natural loss function, demonstrating that energy-based loss functions can produce superior results with significantly fewer unknowns compared to the collocation type methods relying on the residual-based loss functions. Nguyeh-Thanh et al. [35] extended previous work in [34] and applied it to solve nonlinear finite deformation hyperelasticity problems. Additionally, Nguyen-Thanh et al. [36] enhanced convergence by training models in a parametric space, avoiding classical geometry discretization. Using Gauss quadrature in the parametric domain is natural, but mapping between parametric and physical domains poses challenges due to the employment of NURBS basis functions in this approach. Baek et al. [21] proposed a neural network-enhanced reproducing kernel particle method (NN-RKPM) for modeling localization. The RK approximation on a coarse and uniform discretization is employed to approximate the smooth part of the solutions, while the neural networks, determined by the optimization of an energy-based loss function, control parameters to automatically capture the location and oritentaion of the localized solutions. Later on, Baek and Chen [22] propsed an improved version of NN-RKPM in [21] for modeling brittle fracture, in which the NN approximation and background RK approximation are patched together with Partion of Unity to ensure convergence.

In this work, we propose a neural network-enriched Partition of Unity (NN-PU) approximation method for solving PDEs for general elasticity problems, which is a generalization of NN-RKPM previously proposed for modeling localization and fractures [21, 22]. Employing a coarse background discretization, the method utilizes extrinsic NN-based enrichment functions to enhance the background approximation within the Partition of Unity framework [37, 38], in a similar spirit to the hp-Clouds [39, 40], the generalized finite element method [23], and the global-local enrichment [24]. To improve the efficiency and adaptivity of the NN-PU approach, a block-level NN approximation is introduced, where each NN block is designed to target a specific underlying local feature. During offline training, multiple NN enrichment basis sets are pre-trained against different "parent" problems with different local features to embed various underlying solution features into the basis functions. These feature-encoded NN basis sets are then utilized in online simulation stages through transfer learning to capture localized solution efficiently. In this approach, the NN approximation control parameters automatically construct enrichment functions driven by minimizing an energy-based loss function. The proposed NN-PU approximation and feature-encoded transfer learning form an adaptive approximation framework, termed the neural-refinement (n-refinement), departing from traditional polynomial ($p$-based) or mesh ($h$-based) adaptivity complicated by the need to impose conforming constraints.

The remainder of the paper is organized as follows. Section 2 states the problem of interest and outlines the neural network-enriched approximation within the Partition of Unity framework for solving general elasticity problems. Section 3 details NN architectures composed of block-level NN approximation, transferable NN block structure construction, and sub-block architecture forming the NN block. Error analysis of the coupled NN approximation and global NN-PU approximation is also provided in Section 3. Section 4 presents numerical implementations, encompassing offline training of feature-encoded NN enrichment basis sets and the online calculations using trained NN bases via transfer learning. The selection of NN enrichment region based on an energy error indicator is also introduced in Section 4. In Section 5, several numerical examples are presented to assess the solution accuracy, convergence, and effectiveness of the proposed NN-PU approximations. The paper concludes with a discussion and summary in Section 6.

## 2    Neural Network-enriched Partition of Unity Approximation

### 2.1    Model problem

In this study, we consider an elastostatic boundary value problem for easy demonstration of the basic concept. Here a $d$-dimensional elastic body defined in domain $\Omega \in \mathbb{R}^d$ is subjected to boundary conditions on the Dirichlet boundary $\partial\Omega^g$ and Neumann boundary $\partial\Omega^t$ as follows:

$$\boldsymbol{\nabla} \cdot \boldsymbol{\sigma} = \boldsymbol{b} \ \text{ in } \ \Omega$$
$$\boldsymbol{u} = \boldsymbol{g} \ \text{ on } \ \partial\Omega^g \tag{1}$$
$$\boldsymbol{\sigma} \cdot \boldsymbol{n} = \boldsymbol{t} \ \text{ on } \ \partial\Omega^t,$$

where $\boldsymbol{\sigma}$, $\boldsymbol{b}$, $\boldsymbol{g}$, and $\boldsymbol{t}$ denote the Cauchy stress, body force, prescribed displacements, and tractions, respectively. The Cauchy stress for a linear elastic material is defined as:

$$\boldsymbol{\sigma} = 2\mu\boldsymbol{\varepsilon} + \lambda tr(\boldsymbol{\varepsilon})\boldsymbol{I}, \tag{2}$$

with the strain tensor $\boldsymbol{\varepsilon} = \frac{1}{2}(\boldsymbol{\nabla} \otimes \boldsymbol{u} + \boldsymbol{u} \otimes \boldsymbol{\nabla})$, identity tensor $\boldsymbol{I}$, and Lamé's first and second parameters $\lambda$ and $\mu$. The minimization problem corresponds to Eq. (1) is: for $\boldsymbol{u} \in H^1$, $\boldsymbol{u} = \boldsymbol{g}$ on $\partial\Omega^g$,

$$\min_{\boldsymbol{u}} \Pi(\boldsymbol{u}) = \min_{\boldsymbol{u}}[\frac{1}{2}\int_\Omega \boldsymbol{\sigma}\mathbf{:}\boldsymbol{\varepsilon} \, d\Omega - \int_\Omega \boldsymbol{u} \cdot \boldsymbol{b} \, d\Omega - \int_{\partial\Omega^t} \boldsymbol{u} \cdot \boldsymbol{t} \, d\Gamma], \tag{3}$$

where $\Pi(\boldsymbol{u})$ is the potential energy. In this work, $\Pi(\boldsymbol{u})$ will later be used as the loss function for obtaining the neural network enriched solution via loss function minimization.

## 2.2 Partition of Unity approximation with neural network enrichment

Let the domain $\Omega$ be discretized by a set of discrete points $\mathcal{S} \equiv \{1, 2, \cdots, NP\}$, with $NP$ the number of points in the domain discretization, and $V^h = \text{span } \{\Psi_I(\boldsymbol{x})\}_{I \in \mathcal{S}}$ be the finite dimensional approximation space with $\{\Psi_I(\boldsymbol{x})\}_{I \in \mathcal{S}}$ forming a partition of unity (PU), that is, $\sum_{I \in \mathcal{S}} \Psi_I(\boldsymbol{x}) = 1$. Examples of partition of unity functions are the finite element (FE) shape functions, the non-uniform rational basis spline (NURBS) functions in isogeometric analysis (IGA) [41], the moving least-squares (MLS) approximation functions in Element Free Galerkin (EFG) method [42], the reproducing kernel (RK) approximation functions in the reproducing kernel particle method (RKPM) [43, 44], among others. The PU property in the approximation assures first order $L_2$ convergence in solving 2$^\text{nd}$ order PDEs. Higher order convergence can be achieved by enriching the PU via, for example, FE with higher order polynomials, MLS or RK with higher order basis functions as the intrinsic enrichment [45, 46], or both intrinsic and extrinsic enrichments introduced under the frameworks such as the partition of unity method (PUM) [37, 38] and the hp-Clouds [39, 40].

We start with employing the PU functions $\{\Psi_I(\boldsymbol{x})\}_{I\in\mathcal{S}}$ with intrinsic polynomial bases as the "background" approximation, and the neural network (NN) enrichment functions $\{\mathcal{F}_I(\boldsymbol{x})\}_{I\in\mathcal{S}}$ as the extrinsic bases to the background approximation under the PUM [37, 38] framework as follows:

$$u(\boldsymbol{x}) \approx u^h(\boldsymbol{x}) = \sum_{I\in\mathcal{S}} \Psi_I(\boldsymbol{x})\left(\bar{d}_I + \mathcal{F}_I(\boldsymbol{x})\right), \tag{4}$$

where $\bar{d}_I$ is a coefficient associated with $I$-th background PU shape function $\Psi_I(\boldsymbol{x})$. We consider a general neural network approximation represented by:

$$\mathcal{F}_I(\boldsymbol{x}) = \sum_{K=1}^{n_I^\zeta} \zeta_{IK}(\boldsymbol{x})w_{IK} + \beta_I. \tag{5}$$

where $\zeta_{IK}(\boldsymbol{x})$, $n_I^\zeta$, $\beta_I$ are the $K$-th NN basis function, the number of NN basis functions, and the NN bias associated with node $I$, and $w_{IK}$ is the coefficient associated with the NN basis $\zeta_{IK}(\boldsymbol{x})$. Defining a background approximation coefficient $d_I \equiv \bar{d}_I + \beta_I$ to avoid the linear dependency, we express the Neural Network-enriched Partition of Unity (NN-PU) approximation as:

$$u^h(\boldsymbol{x}) = \sum_{I\in\mathcal{S}} \Psi_I(\boldsymbol{x})\left(d_I + \sum_{K=1}^{n_I^\zeta} \zeta_{IK}(\boldsymbol{x})w_{IK}\right). \tag{6}$$

In this work, we call $\zeta_{IK}(\boldsymbol{x})$ the *NN enrichment basis*, and $w_{IK}$ the *NN enrichment coefficient*.

The approximation $u^h(\boldsymbol{x})$ can also be viewed as the superposition of a background approximation $u^{BG}(\boldsymbol{x})$ and an NN approximation $u^{NN}(\boldsymbol{x})$, and $\Psi_I(\boldsymbol{x})$ with the partition of unity property serves as a patching function for the NN enrichment bases $\{\zeta_{IK}(\boldsymbol{x})\}_{K=1}^{n_I^\zeta}$. Rewrite Eq. (6) as:

$$u^h(\boldsymbol{x}) = u^{BG} + u^{NN}, \tag{7}$$

with

$$u^{BG}(\boldsymbol{x}) = \sum_{I \in \mathcal{S}} \Psi_I(\boldsymbol{x}) d_I,$$

$$u^{NN}(\boldsymbol{x}) = \sum_{I \in \mathcal{S}^\zeta} \Psi_I(\boldsymbol{x}) \sum_{K=1}^{n_I^\zeta} \zeta_{IK}(\boldsymbol{x}) w_{IK}, \tag{8}$$

where $\mathcal{S}^\zeta$ is the set of nodes with NN enrichment, $\mathcal{S}^\zeta \subset \mathcal{S}$, and $w_{IK} = 0 \ \forall \ I \in \mathcal{S} \backslash \mathcal{S}^\zeta$, to allow for enrichment only at selective set of points, which will be discussed in Sec. 4.1.3. Let's now consider global NN-enrichment bases that are not associated with the background discretization, i.e., $\zeta_K(\boldsymbol{x}) \equiv \zeta_{IK}(\boldsymbol{x})$ and $n^\zeta \equiv n_I^\zeta$, independent of the discretization index $I$. Then, $u^{NN}(\boldsymbol{x})$ defined in Eq. (8) can be written as:

$$u^{NN}(\boldsymbol{x}) = \sum_{I \in \mathcal{S}^\zeta} \sum_{K=1}^{n^\zeta} \zeta_K(\boldsymbol{x}) \Psi_I(\boldsymbol{x}) w_{IK}, \tag{9}$$

where $n^\zeta$ is the number of NN enrichment bases and let $\mathbf{W}_K = \{w_{IK}\}_{I \in \mathcal{S}}$ be an NN enrichment coefficient set.

Note that the proposed approximation takes a similar form to the standard Partition of Unity and $hp$-clouds methods [37-40]. A key difference being in its flexibility in dynamically adapting enrichment functions constructed via a neural network loss function minimization to be discussed in the next section without pre-defined enrichment functions. We term this type of refinement a neural-refinement ($n$-refinement). The proposed approach is also fundamentally different from solving PDEs by the Galerkin approximation at a stationary, which instead is solved by a Ritz-type approach through a neural-network-based potential function minimization to be discussed in the next section.

***Remark 2.1.*** While the background shape functions $\Psi_I(\boldsymbol{x})$ can be of any type that possesses the partition of unity, this work adopts the reproducing kernel (RK) functions [43, 44] for the background approximations. RK approximation offers independent choice over the order of continuity (smoothness) and the order of completeness, enabling high-order continuity with lower-order bases, particularly beneficial for approximating smooth solutions while replying on the NN approximations to capture localized behavior such as corner singularities, see [43, 44, 47, 48] for details.

***Remark 2.2.*** The computational architecture of the neural network approximation in Eq. (9) can be decomposed into multiple NN blocks as those designed in [21, 22]. On this concept of block-level NN approximation, we can rewrite Eq. (9) with $n^B$ blocks as follows:

$$u^{NN}(\boldsymbol{x}) = \sum_{J=1}^{n^B} u_J^B(\boldsymbol{x}) = \sum_{J=1}^{n^B} \left( \sum_{I \in \mathcal{S}_J^\zeta} \sum_{L=1}^{n_J^\zeta} \zeta_{JL}(\boldsymbol{x}) \Psi_I(\boldsymbol{x}) w_{IJL} \right), \qquad (10)$$

where $u_J^B(\boldsymbol{x})$, $n_J^\zeta$, $\zeta_{JL}$, and $\mathcal{S}_J^\zeta$ denote the block-level NN approximation, the number of NN bases, $L$-th NN basis, and the NN enriched nodeset of NN block $J$, respectively.

# 3 Transferable Neural Network Architecture for NN Approximation

## 3.1 Block-level NN approximation

One goal of this work is to develop feature-encoded and computationally efficient NN enrichment bases based on the neural network block architectures. The block-level NN approximation, as shown in Figure 1(b), possesses high sparsity in the network structure compared to the densely connected counterpart shown in Figure 1(a). Meanwhile, a block neural network has the flexibility to incorporate multiple feature-encoded NN bases targeting different localized features in the solution. This is a unique feature for constructing local enrichment for enhanced solution accuracy without the conventional mesh adaptive refinement.

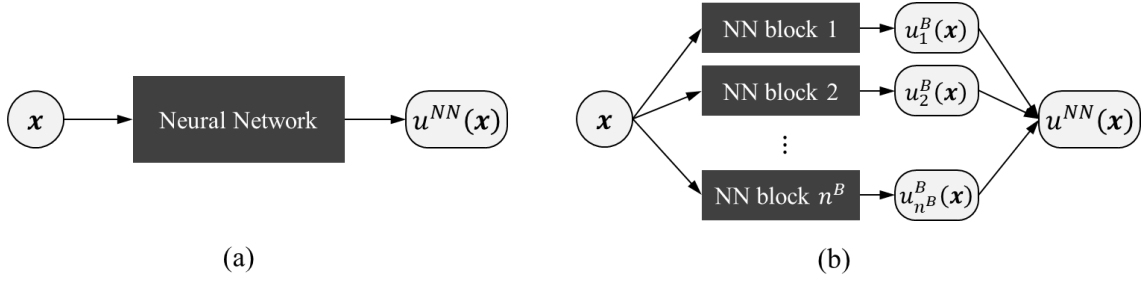

**Figure 1: Neural network architectures for NN approximation: (a) a single-block densely connected deep neural network and (b) multi-block sparsely connected neural network**

## 3.2 Transferable neural network block architecture

Concerning the network architecture of each neural network block, it is not feasible to optimize the entire block structure during the online computation due to the associated high computational cost. Conversely, utilizing fully trained (fixed) neural network block is also undesirable, as they may lack the flexibility in adjusting enrichment functions to capture unseen solution behaviors in the online computation. We aim to construct partially trained NN block structures. Herein, partial training entails training enrichment functions in an offline stage and employing them for online computations. This transfer learning-type approach aims to achieve a balance between computational accuracy and efficiency.

To this end, we define a block-level NN approximation as the composition of offline-trained NN sub-block $\mathcal{N}^\zeta$ and online-optimized NN sub-blocks $\mathcal{N}^P$ and $\mathcal{N}^C$ as follows (omitting the subscript $J$ associated with the $J$-th NN block for brevity):

$$u^B(\boldsymbol{x}) = \mathcal{N}^C(\cdot; \mathbf{W}^C) \circ \mathcal{N}^\zeta(\cdot; \mathbf{W}^\zeta) \circ \mathcal{N}^P(\boldsymbol{x}; \mathbf{W}^P), \tag{11}$$

where $\circ$ denotes the function composition operator, $\mathcal{N}^P$, $\mathcal{N}^\zeta$, and $\mathcal{N}^C$ denote the *Parametric Sub-block*, *NN Basis Sub-block*, and *Coefficient Sub-block*, as shown in Figure 2, and $\mathbf{W}^P$, $\mathbf{W}^\zeta$, and $\mathbf{W}^C$ are the weight sets associated with $\mathcal{N}^P$, $\mathcal{N}^\zeta$, and $\mathcal{N}^C$, respectively.

The Parametric Sub-block, $\mathcal{N}^P(\cdot, \mathbf{W}^P): \boldsymbol{x} \to \boldsymbol{y}$, is designed to generate parametric coordinates $\boldsymbol{y}$ to facilitate the utilization of pre-trained NN bases in a parametric coordinate $\boldsymbol{y}$ to represent functions with local features. Nonlinear mapping $\boldsymbol{x} \to \boldsymbol{y}$ also allows a representation of complex high-dimensional features by NN bases in a lower-dimensional manifold. The NN Basis Sub-block, $\mathcal{N}^\zeta(\cdot; \mathbf{W}^\zeta): \mathbf{y} \to \mathbf{Z}$, takes the parametric coordinates $\mathbf{y} = [y_1, \cdots, y_d]$ and generates a NN basis vector $\mathbf{Z} = [\zeta_1, \cdots, \zeta_{n^\zeta}]$. The weight set $\mathbf{W}^\zeta$ is trained in the offline stage and can be further optimized during the online calculation. This sub-block is designed to embed the trained NN bases with underlying local features that are not captured by the background approximation functions. The linear combination of the NN basis by the coefficient set $\mathbf{W}^C$ is also computed during the online stage. With the computed $\mathbf{W}^C$, the Coefficient Sub-block, $\mathcal{N}^C(\cdot, \mathbf{W}^C): \mathbf{Z} \to u^B$, adds the contribution of each NN enrichment basis to the total solution of the problem.

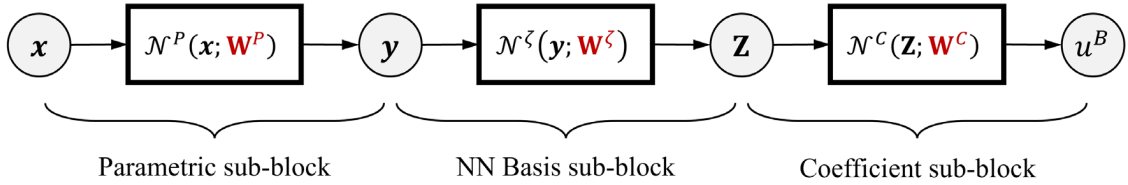

**Figure 2: Global architecture of a transferable neural network adopted for a single NN block**

### 3.3 Offline training of single feature by a reduced NN block architecture

In order to capture various types of solution features, an approximation may be constructed upon multiple NN blocks during online computation, each of which is partially trained against a specific solution feature in the offline training stage. In this work, a "parent" problem that presents a distinct solution feature to embed is used to train the NN Basis Sub-block $\mathcal{N}^\zeta$ with a single NN block during the offline training stage, and the physical domain of the "parent" problem is utilized as the parametric coordinates. During the offline training, the Parametric Sub-block is set to be the identity map, i.e., $\mathcal{N}^P: \mathbf{x} \to \mathbf{x}$. Then, the block-level approximation in Eq. (11) can be written as:

$$u^{NN}(\mathbf{x}) = u^B(\mathbf{x}) = \mathcal{N}^C(\cdot; \mathbf{W}^C) \circ \mathcal{N}^\zeta(\mathbf{x}; \mathbf{W}^{\zeta^0}) \equiv \mathcal{N}^C(\mathbf{Z}^0; \mathbf{W}^C), \quad (12)$$

and the *parent* NN basis set $\mathbf{Z}^0 = (\zeta_1^0, \cdots, \zeta_{n^\zeta}^0)$ is defined as:

$$\mathbf{Z}^0 = \mathcal{N}^\zeta\big(\boldsymbol{x}; \mathbf{W}^{\zeta^0}\big), \tag{13}$$

where the NN Basis Sub-block $\mathcal{N}^\zeta$ directly takes the physical coordinates of the "parent" problem with single feature as its input and $\mathbf{W}^{\zeta^0}$ is the converged weight sets of NN Basis Sub-block during the offline training. Figure 3 depicts the reduced NN architecture for an offline training with single feature. Note that setting $\mathcal{N}^P$ to the identity map eliminates the operations associated with $\mathcal{N}^P$, which reduces the computational cost during the offline stage. Meanwhile, although the Coefficient Sub-block $\mathcal{N}^C$ is involved, and its weight set $\mathbf{W}^C$ is optimized together with $\mathbf{W}^\zeta$ in the offline training, $\mathbf{W}^C$ obtained in the offline training is discarded and re-evaluated in the online simulation as the online computation could involve multiple features, and the weights for the linear combination of NN bases need to be recomputed.

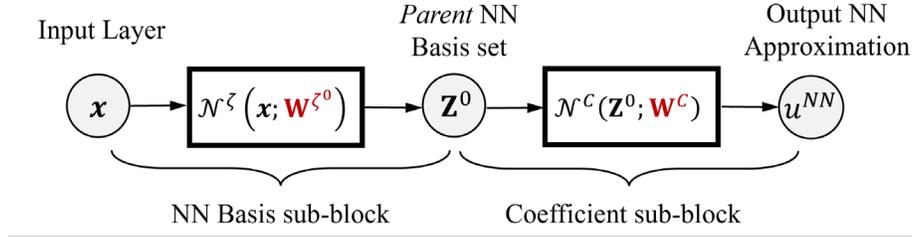

**Figure 3: Reduced NN architecture for one offline training**

### 3.4 Structures of NN sub-block architectures

#### 3.4.1 Architecture of $\mathcal{N}^\zeta$

The NN basis Sub-block is constructed with the standard deep neural network (DNN) architecture. Figure 4 shows a general DNN architecture, represented by $\mathcal{N}^{DNN}(\boldsymbol{\mathcal{I}}; \boldsymbol{w})$ with input vector $\boldsymbol{\mathcal{I}} = \big(i_1, \cdots, i_{n_i}\big)$, output vector $\boldsymbol{\mathcal{O}} = \big(o_1, \cdots, o_{n_o}\big)$, and weight set $\boldsymbol{w}$. The weight set $\boldsymbol{w}$ is defined as the collection of the weights and biases of all the layers in DNN, i.e., $\boldsymbol{w} = \{\boldsymbol{\theta}_\ell, \beta_\ell\}_{\ell=1}^{NL}$, where $\boldsymbol{\theta}_\ell$, $\beta_\ell$, and $NL$ are the weight matrix of layer $\ell$, the bias of layer $\ell$, and the total number of hidden and output layers, respectively. Note that layer 0 is the input layer, layer $NL$ is the output layer, and the intermediate layers are hidden layers of DNN. The intermediate output of layer $\ell$, denoted as $\boldsymbol{h}_\ell$, is written as:

$$\boldsymbol{h}_\ell = a_\ell(\boldsymbol{\mathcal{L}}(\boldsymbol{h}_{\ell-1}; \{\boldsymbol{\theta}_\ell, \beta_\ell\})), \ \ell = 1, \cdots, NL - 1, \tag{14}$$

where the linear combination operator $\mathcal{L}(\boldsymbol{h}; \{\boldsymbol{\theta}, \beta\})$ is defined as:

$$\mathcal{L}(\boldsymbol{h}; \{\boldsymbol{\theta}, \beta\}) = \boldsymbol{\theta}\boldsymbol{h} + \beta, \tag{15}$$

and $a_\ell$ denotes the activation function applied to layer $\ell$, with exponential linear unit (ELU) activation function [49] considered in this work. Note that we have $\boldsymbol{\mathcal{I}} = \boldsymbol{h}_0$ and $\boldsymbol{\mathcal{O}} = \boldsymbol{h}_{NL}$. For the NN Basis Sub-block $\mathcal{N}^\zeta$, a DNN is adopted with input $\boldsymbol{\mathcal{I}} = \boldsymbol{y}$, with $\boldsymbol{y}$ the parametric coordinates, output $\boldsymbol{\mathcal{O}} = \mathbf{Z}$, with $\mathbf{Z}$ the NN basis set, and weight set $\boldsymbol{w} = \mathbf{W}^\zeta$.

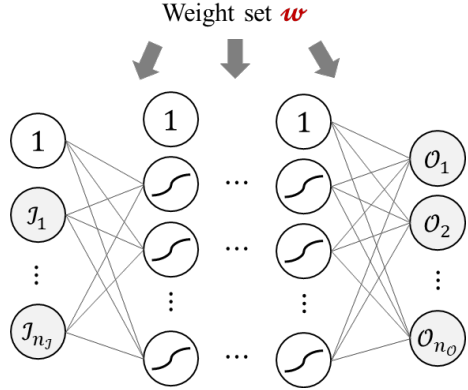

**Figure 4: A NN basis Sub-block $\mathcal{N}^\zeta$ constructed by a densely connected neural network**

### 3.4.2 Architecture of $\mathcal{N}^P$

As discussed in Section 3.3, the Parametric Sub-block $\mathcal{N}^P$ is set to be the identity map for the offline training stage. In this work, $\mathcal{N}^P$ is inspired by the residual neural network (ResNet) [50]:

$$\mathcal{N}^P(\boldsymbol{x}; \mathbf{W}^P) = \mathcal{N}^{DNN}(\boldsymbol{x}; \mathbf{W}^P) + \boldsymbol{x}, \tag{16}$$

and the network architecture is shown in Figure 5. Note that by setting $\mathbf{W}^P = \mathbf{0}$, the identity map used for the offline training stage is recovered. Besides its identity reproducibility, ResNet is proven to avoid the gradient vanishing issues commonly encountered in deep neural networks, thereby enhancing computational efficiency [50].

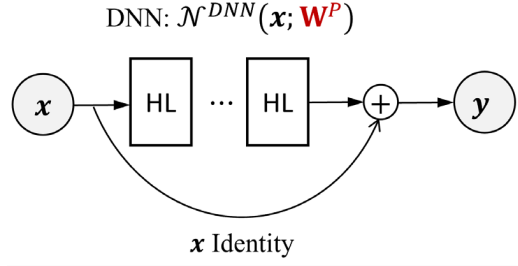

**Figure 5: Modified architecture for Parametrization Network** (HL stands for Hidden Layer)

### 3.4.3 Architecture of $\mathcal{N}^C$

The Coefficient Sub-block $\mathcal{N}^C$ is defined as follows:

$$\mathcal{N}^C(\mathbf{Z}; \mathbf{W}^C) = \sum_{I \in \mathcal{S}^\zeta} \Psi_I(\boldsymbol{x}) \mathcal{F}_I(\boldsymbol{x}; \mathbf{W}^C), \tag{17}$$

with $\mathbf{Z} = [\zeta_1, \cdots, \zeta_{n^\zeta}]$ and $\mathcal{F}_I(\boldsymbol{x}; \mathbf{W}^C) = \sum_{L=1}^{n^\zeta} \zeta_L(\boldsymbol{x}) w_{IL}$. The sub-block architecture is shown in Figure 6.

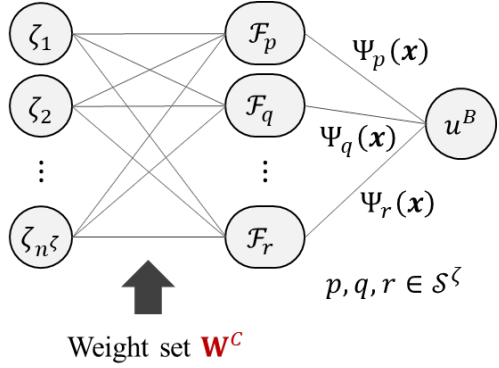

Weight set $\mathbf{W}^C$

**Figure 6: Network architecture of Output Sub-block**

## 3.5 Convergence

Melenk and Babuska (1996) [37] provided the following global error bound of a Partition of Unity approximation as follows: with a generic constant $C_G$,

$$\|u^h - u\|_{0,\Omega} \leq C_G \left( \sum_I \epsilon_I^2 \right)^{1/2}, \tag{18}$$

where $\epsilon_I$ is a local error limit and bounds the error on $\Omega_I \equiv \text{supp}(\Psi_I(x))$:

$$\|w_I - u\|_{0,\Omega_I} \leq \epsilon_I, \tag{19}$$

where $w_I = w_I(x)$ is the local approximation function of a general Partition of Unity approximation, i.e., $u^h(x) = \sum_I \Psi_I(x) w_I(x)$. Let $p^\Psi$ and $p^w$ be the order of complete polynomial reproduced by $\text{span}\{\Psi_I\}_{I \in \mathcal{S}}$ and $w_I(x)$, respectively. For the case that $p^w$ is the leading order, i.e., $p^w \geq p^\Psi$, Duarte and Oden (1996) [40] proved the following global error estimate for smooth $u$:

$$\|u^h - u\|_{0,\Omega} \leq \mathcal{O}(h^{p^w+1}), \tag{20}$$

where $h = \max\limits_I h_I$ is the maximum nodal spacing. For the case that $p^\Psi$ is the leading order and $w_I(x) = d_I$ is constant, e.g., the Reproducing Kernel approximation, Han et al. (2002) [51] and Chen et al. (2003) [52] showed the following global error estimate for smooth $u$:

$$\|u^h - u\|_{0,\Omega} \leq \mathcal{O}(h^{p^\Psi+1}). \tag{21}$$

### 3.5.1 Convergence of pure neural network approximation in Partition of Unity

Based on the universal approximation theorem [2] and a more recent study in [53], the local error of a neural network enrichment function $\mathcal{F}_I = \mathcal{F}(x; \mathbf{W}_I)$ can be represented as follows: for $u \in L_2$,

$$\|\mathcal{F}_I - u\|_{0,\Omega_I} \leq C_{\Omega_I} \hat{n}^{1/2}, \tag{22}$$

where $\hat{n}$ is a size measure of the neural network and $C_{\Omega_I}$ is a parameter defined as

$$C_{\Omega_I} = \left( \int_{\Omega_I} B|\hat{u}(\omega)| \, d\Omega(\omega) \right)^{1/2}, \tag{23}$$

where $\hat{u}$ is the frequency-domain solution, $\omega$ is the frequency, and $B$ is a constant independent of the network size measure $\hat{n}$ and the nodal spacing of $I$-th background node $h_I$. To remove the dependency of $C_{\Omega_I}$ on the support size of $\Psi_I$, the physical coordinate can be normalized by defining a parametric coordinate $\xi \equiv (x - x_I)/h_I$ with parametric domain $\boxdot$, which yields:

$$C_{\Omega_I} = h_I^{d/2} \left( \int_{\boxdot} B|\hat{u}(\omega)| \, d\xi(\omega) \right)^{1/2} \equiv h_I^{d/2} \hat{C}, \tag{24}$$

where $d$ is the space dimension. With Eqs. (18), (22), and (24), we have the following relation:

$$\|u^{NN} - u\|_{0,\Omega} \leq C_G \hat{C} \left( \sum_I h_I^d \right)^{1/2} \hat{n}^{1/2} \leq C_G C \left( (NP) h^d \right)^{1/2} \hat{n}^{1/2} \tag{25}$$

With $NP \sim \mathcal{O}(h^{-d})$, Eq. (25) becomes:

$$\|u^{NN} - u\|_{0,\Omega} \leq \mathcal{O}(\hat{n}^{1/2}). \tag{26}$$

### 3.5.2 Convergence of NN-PU approximation

Consider an NN-PU approximation with $u^{BG}$ of degree $p$ and $u^{NN}$ of size measure $\hat{n}$. In order to estimate the local error bound of NN-PU, a local interpolation error from $u^{BG}$ is first estimated, followed by the estimation of the error from $u^{NN}$. Note that in this work $u^{BG} = u^{RK}$. Let $L^p(x; \mathbf{C})$ be a degree-$p$ polynomial function of $x$ with a set of constants $\mathbf{C}$. For $u^{BG}$ that can locally reproduce a complete polynomial of degree $p$, and considering the Taylor expansion of $u$, we have the local approximation:

$$
\begin{aligned}
&\|u^{BG} + u^{NN} - u\|_{0,\Omega_I} \\
&= \|L^p(x; \mathbf{C}_I) + \mathcal{F}(x; \mathbf{W}_I) - u(x)\|_{0,\Omega_I} \\
&= \left\| \mathcal{F}(x; \mathbf{W}_I) - \frac{(x - \bar{x})^{p+1}}{(p+1)!} u^{(p+1)}(\bar{x}) - \mathcal{O}((x - \bar{x})^{p+2}) \right\|_{0,\Omega_I}.
\end{aligned}
\tag{27}
$$

Further introducing a parametric coordinate $\xi \equiv (x - \bar{x})/h_I$ with parametric domain $\boxdot$, we have:

$$\left\| \mathcal{F}(x; \mathbf{W}_I) - \frac{(x - \bar{x})^{p+1}}{(p+1)!} u^{(p+1)}(\bar{x}) - \mathcal{O}((x - \bar{x})^{p+2}) \right\|_{0,\Omega_I}$$

$$= h_I^{d/2} \left\| \mathcal{F}(\xi; \overline{\mathbf{W}}_I) - h_I^{p+1} \frac{\xi^{p+1}}{(p+1)!} u^{(p+1)}(\bar{x}) - \mathcal{O}(h_I^{p+2} \xi^{p+2}) \right\|_{0,\boxdot} \tag{28}$$

$$= h_I^{p+1+d/2} \left\| \mathcal{F}(\xi; \overline{\overline{\mathbf{W}}}_I) - \frac{\xi^{p+1}}{(p+1)!} u^{(p+1)}(\bar{x}) - \mathcal{O}(h_I \xi^{p+2}) \right\|_{0,\boxdot}.$$

Utilizing Eq. (22), we obtain the following local error estimate when $h_I \to 0$: for $u^{(p+1)} \in L_2$,

$$\| u^{BG} + u^{NN} - u \|_{0,\Omega_I} \leq h_I^{p+1+d/2} \left\| \mathcal{F}(\xi; \overline{\overline{\mathbf{W}}}_I) - \frac{\xi^{p+1}}{(p+1)!} u^{(p+1)}(\bar{x}) \right\|_{0,\boxdot} \tag{29}$$

$$\leq \hat{C}_I h_I^{p+1+d/2} \hat{n}^{1/2} \equiv \epsilon_I^{NNPU}.$$

With Eq. (18), Eq. (29), and $NP \sim \mathcal{O}(h^{-d})$, we have the following global NN-PU approximation error with background RK approximation:

$$\| u^h - u \|_{0,\Omega} \leq \mathcal{O}(h_I^{p+1} \hat{n}^{1/2}). \tag{30}$$

## 4    Numerical Implementation

The minimization problem in Eq. (3) can be expressed as:

$$\min_{\mathbf{d}, \mathbb{W}} \left[ \Pi\left( \boldsymbol{u}^h(\mathbf{d}, \mathbb{W}) \right) \right], \tag{31}$$

where $\boldsymbol{u}^h(\mathbf{d}, \mathbb{W}) = \boldsymbol{u}^{RK}(\mathbf{d}) + \boldsymbol{u}^{NN}(\mathbb{W})$ is the NN-PU approximation with background RK approximation $\boldsymbol{u}^{RK}$ and $\mathbb{W} = \{\mathbb{W}^P, \mathbb{W}^\zeta, \mathbb{W}^C\}$ denoting the neural network weight sets of Parametric, NN Basis, and Coefficient sub-blocks, respectively, with $\mathbb{W}^P = \{\mathbf{W}_J^P\}_{J=1}^{n^B}$, $\mathbb{W}^\zeta = \{\mathbf{W}_J^\zeta\}_{J=1}^{n^B}$, and $\mathbb{W}^C = \{\mathbf{W}_J^C\}_{J=1}^{n^B}$. This work follows a staggered approach to obtain the NN-PU approximations, where the background RK solutions are first computed without NN enrichments, followed by the optimization of NN-PU approximations with fixed pre-computed background solutions. Figure 7 outlines the solution procedures for obtaining the NN-PU approximation.

## 4.1 NN-PU approximation by loss function minimization

### 4.1.1 Input variables

The required inputs for the NN-PU optimization are given in the gray block in Figure 7. The background RK coefficient vector $\mathbf{d} \equiv \{d_I\}_{I \in \mathcal{S}}$ is obtained by the standard Galerkin approximation with RK approximation functions [43, 44, 54]. Hence, the physical coordinates of domain integration evaluation points $\{x_e\}_{e=1}^{Ne}$ and a set of RK shape functions evaluated at those evaluation points $\{\{\Psi_I(x_e)\}_{I=1}^{NP}\}_{e=1}^{Ne}$ are given as inputs. The background RK solutions are also pre-evaluated at those evaluation points and inputted for the computation of NN-PU approximation. Background nodal coordinates and RK nodal shape functions, evaluated with normalized support sizes of $a$ and $2a$, are also pre-determined for error estimation analysis and updating NN-enriched nodesets as detailed in Sec. 4.1.3.

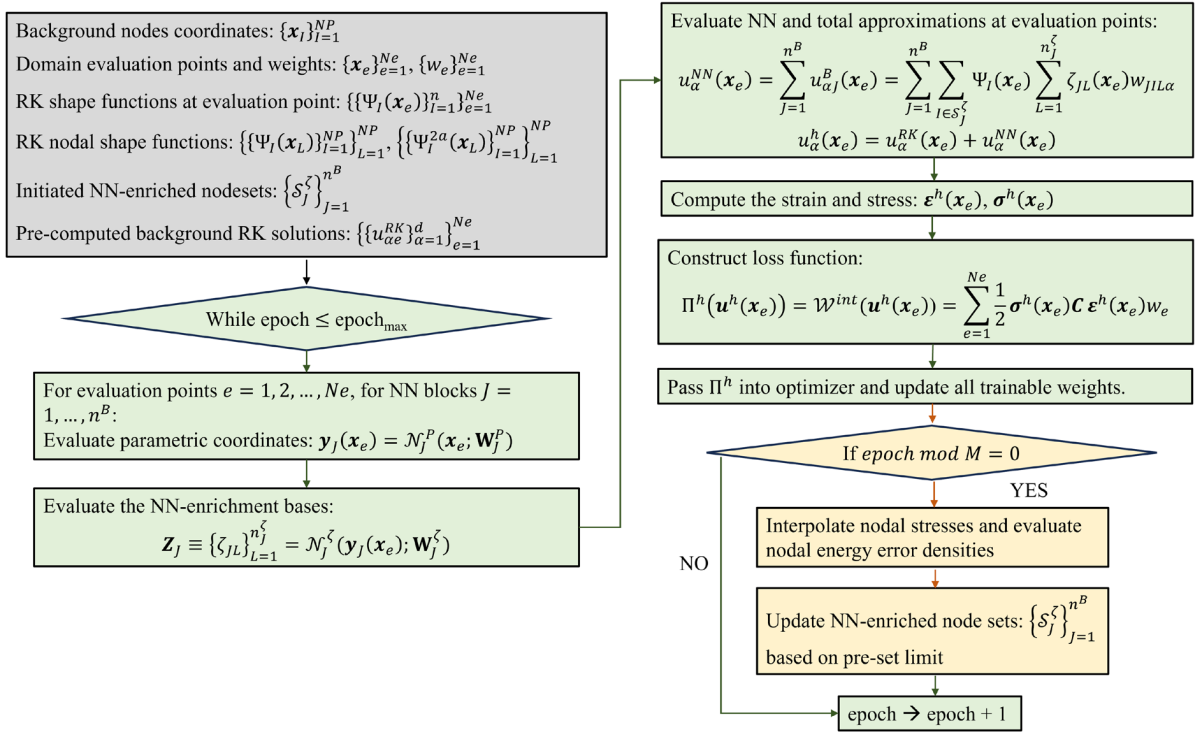

**Figure 7: Flowchart of the NN-PU solution procedures** (gray block for inputs, green blocks for loss function minimization procedures, yellow blocks for updating NN-enriched nodesets)

### 4.1.2 Loss function minimization

The NN loss function minimization procedures are outlined in green blocks in Figure 7. For the *offline training* of the "parent" NN enrichment basis functions given background RK solution $\boldsymbol{u}^{RK}$, the weights and biases of the NN Basis and Coefficient sub-Blocks are obtained by:

$$\mathbb{W}^{\zeta^0}, \mathbb{W}^C = \arg\min_{\mathbb{W}^{\zeta^0}, \mathbb{W}^C} \left[ \Pi(\boldsymbol{u}^{RK}(\mathbf{d}) + \boldsymbol{u}^{NN}(\{\mathbb{W}^{\zeta^0}, \mathbb{W}^C\})) \right] \tag{32}$$

subjected to $\boldsymbol{u}^{RK}(\boldsymbol{x}) = \boldsymbol{g}(\boldsymbol{x})$, $\boldsymbol{u}^{NN}(\boldsymbol{x}) = \boldsymbol{0}$, $\forall \boldsymbol{x} \in \partial\Omega^g$.

Recall that the Parametric sub-Block is designated as an identity map during offline training so that $\mathbb{W}^P = \boldsymbol{0}$. The weight set $\mathbb{W}^{\zeta^0}$ of the NN Basis sub-Block and the resultant "parent" NN enrichment basis sets $\mathbb{Z}^0 = \left\{ \mathcal{N}_J^\zeta \left( \boldsymbol{x}; \mathbb{W}_J^{\zeta^0} \right) \right\}_{J=1}^{n^B}$ are preserved for subsequent online computations for problems with similar local features. As indicated in Eq. (32), the essential boundary conditions are satisfied by the background RK solutions. The homogenous essential boundary condition is imposed on NN approximations by setting $w_{IJL} = 0$ for all $I$ on essential boundaries.

The *online computation* of problems involving multiple features employ transfer-learning with feature-encoded "parent" NN enrichment basis sets, each targeting a single feature, as follows:

*Given $\boldsymbol{u}^{RK}(\boldsymbol{d})$ and $\mathbb{Z}^0 = \left\{ \mathcal{N}_J^{\zeta} \left( \boldsymbol{x}; \boldsymbol{W}_J^{\zeta^0} \right) \right\}_{J=1}^{n^B}$, find $\{\mathbb{W}^P, \mathbb{W}^C\}$ such that the loss function $\Pi$ is minimized as:*

$$\{\mathbb{W}^P, \mathbb{W}^C\} = \arg\min_{\{\mathbb{W}^P, \mathbb{W}^C\}} \left[ \Pi(\boldsymbol{u}^{RK}(\boldsymbol{d}) + \boldsymbol{u}^{NN}(\{\mathbb{W}^P, \mathbb{W}^{\zeta^0}, \mathbb{W}^C\})) \right]$$

(33)

*subjected to $\boldsymbol{u}^{RK}(\boldsymbol{x}) = \boldsymbol{g}(\boldsymbol{x})$, $\boldsymbol{u}^{NN}(\boldsymbol{x}) = \boldsymbol{0}$, $\forall \boldsymbol{x} \in \partial\Omega^g$.*

During online solution computations, weight sets of NN Basis sub-Blocks are known and pre-trained in a parametric space during offline computation, while the to-be-optimized weight set $\mathbb{W}^P$ finds an optimal mapping between $\boldsymbol{x}$ and $\boldsymbol{y}_J$ given the NN enrichment basis sets obtained from "parent" problems, and $\mathbb{W}^C$ is responsible for the coefficients of the NN basis functions.

### 4.1.3 Error indicator based on Reproducing Kernel Background Approximation

In the online computation, the NN enrichment nodes are initially informed by the error indicator for RK background solution. It has been investigated in [55, 56], the discrete convolution of two reproducing kernels with $n$th-order bases and support sizes $a$ results in a filtered kernel of the same order of completeness. This property makes reproducing kernel discrete convolution suitable for indicating errors and identifying high gradients in the RK background approximation needing solution enrichments. The filtered RK approximated stress, or recovered stress [56], denoted as $\boldsymbol{\sigma}^*(\boldsymbol{x})$, is obtained as follows:

$$\boldsymbol{\sigma}^*(\boldsymbol{x}) = \sum_{I \in \mathcal{S}} \Psi_I^{2a}(\boldsymbol{x}) \boldsymbol{\sigma}^a(\boldsymbol{x}_I),$$

(34)

where $\boldsymbol{\sigma}^a(\boldsymbol{x}_I)$ is the nodal RK stress solution obtained by RK kernel function with support "a", and $\{\Psi_I^{2a}\}_{I\in\mathcal{S}}$ is the set of RK shape functions constructed with twice the support size used for the RK solution. The error indication of the stress solution is defined as:

$$\boldsymbol{e}_{\boldsymbol{\sigma}}^* = \boldsymbol{\sigma}^* - \boldsymbol{\sigma}^a, \tag{35}$$

and the nodal energy error density $\rho_I^*$ can be computed using $\boldsymbol{e}_{\boldsymbol{\sigma}}^*$ as follows:

$$\rho_I^* = \left\{\frac{1}{2}\boldsymbol{e}_{\boldsymbol{\sigma}}^*(\boldsymbol{x}_I)^T \boldsymbol{C}^{-1}\boldsymbol{e}_{\boldsymbol{\sigma}}^*(\boldsymbol{x}_I)\right\}^{1/2}, \qquad \forall I \in \mathcal{S}, \tag{36}$$

where $\boldsymbol{C}$ is the elasticity tensor. The NN enriched nodeset $S^\zeta$ for a NN block is defined based on the nodal energy error densities:

$$S^\zeta = \{I | \rho_I^* \geq \tau \cdot \rho_{max}^*\}, \tag{37}$$

where $\rho_{max}^* \equiv \max_{I\in\mathcal{S}} \rho_I^*$, and $\boldsymbol{\tau} \in (0,1]$ is a scaling parameter. Note that nodes situated on the Dirichlet boundaries are excluded in the node set $S^\zeta$ for NN enrichment.

As per Eq. (9), the size of the NN enrichment coefficient set is closely related to the size of the NN enriched nodeset. Further, the error estimation analysis outlined in Eqs. (34)-(36) is performed periodically throughout the loss function minimization process, outlined in the yellow block in Figure 7. The NN-enriched nodesets $S^\zeta$ determined by Eq. (37) are updated in the gradient decent iteration until they stabilize.

# 5   Numerical Examples

In this section, we present the offline construction of NN enrichment functions within the "parent" problems and illustrate their utilization in online solution computations through transfer learning. Subsequent subsections showcase the numerical procedures and explore solution convergence when employing the suggested NN enrichment for problems involving localized features:

- Section 5.1: Offline training of NN enrichment functions for "parent" problems with local features

- Section 5.2: Online solution of problems with similar local features by transfer leaning

- Section 5.3: Online solution of problems with multiple local features by transfer leaning

Unless specified otherwise, the background RK shape functions are constructed using the linear basis and cubic B-spline kernel function with a normalized support size (normalized with nodal distance) $a = 2.0$. Stabilized Conforming Nodal Integration [54] is employed for domain integration of the background solution. The RK shape function scaling approach [57] is introduced for strong imposition of essential boundary conditions on the pre-computed background RK approximations. Zero NN coefficients are enforced on essential boundaries to ensure that the NN solutions satisfy the homogenous essential boundary condition. Material properties $E = 10^4$ and $\nu = 0.3$ are chosen for all examples in this section. In the NN Basis sub-block, the exponential linear unit (ELU) activation function [49] is employed for non-linear transformations, effectively addressing the vanishing gradient issue commonly associated with ReLU activation functions. ELU's distinctive inflection point in its first derivative, as shown in Figure 8, makes it suitable for modeling stress singularity or concentration. The hyperbolic tangent activation function is selected for all hidden layers in the Parametric sub-block except the last one, where a linear activation function is applied. NN optimization utilizes the first-order Adam (adaptive momentum) optimizer [58], provided in the open-source machine learning library TensorFlow [26].

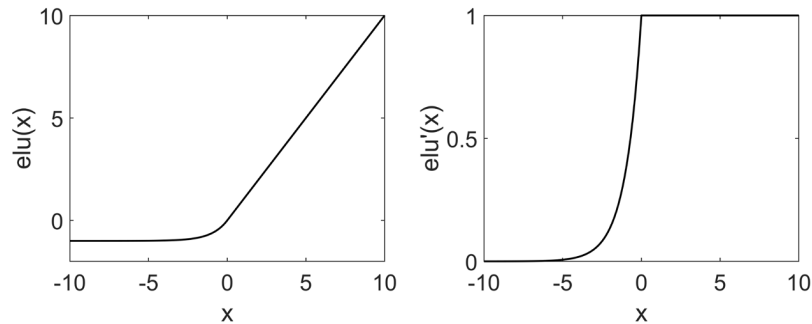

**Figure 8: Plot of the ELU activation function (left) and its first derivative (right)**

In the following demonstration problems, the solutions obtained from NN-PU with NN-enrichment will be compared with finite element (FE) solutions with various levels of *h-* and *p-* refinements. *Here, the FEM solutions are obtained by Galerkin approximation* *of static equilibrium, that is the stationary of potential energy, while the NN-PU solution are obtained via minimization of potential energy through gradient decent iterative procedures as stated in Section 4.1.* In addition, according to the property that the Galerkin approximation (obtained from the stationary state of potential energy) overestimates the potential energy of the continuum system in equilibrium [59], the potential energy in the NN-PU loss function will be compared with the FE potential energy as an accuracy measure. The NN-PU CPU reported in the numerical examples is with both background RK solution computations and online NN optimizations included.

## 5.1    Offline training of parent NN enrichment functions

This subsection demonstrates the offline training of NN enrichment functions in two "parent" problems, and examines the efficiency, accuracy, and convergence of the proposed NN-PU computational framework. These problems target local features and employ a single NN block with the Parametric sub-Block set to be an identity map in constructing the NN enrichment functions. Neural network iteration convergence is met either when the loss function reduction is less than 10 percent of the learning rate in 5 epochs or when the total epochs surpass 30,000.

### 5.1.1    Parent Problem 1: Plate with a hole

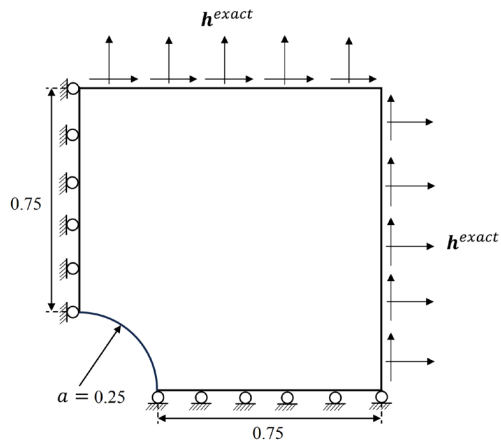

**Figure 9: Plate with a hole unit problem setting**

The first "parent" problem concerns a plate with a circular hole subjected to a far-field uniaxial tension of 10. Due to the axisymmetry, only the upper right quadrant of the domain is modeled, as shown in Figure 9. Symmetric boundary conditions are imposed on the left vertical and bottom horizontal edges, while the exact traction field, calculated using the analytical solutions presented in [60], is prescribed on the remaining edges. This problem focuses on analyzing the convergence rates of NN-PU solutions and how the selection of the enrichment region influences the potential's minimizing performance.

To facilitate the convergence studies, three levels of background discretization are employed, as shown in Figure 10, and a single hidden layer with 5 NN enrichment basis functions is used for constructing the NN basis functions. The solution accuracy is investigated in terms of the normalized displacement error norms ($e_{L2}$) as follows:

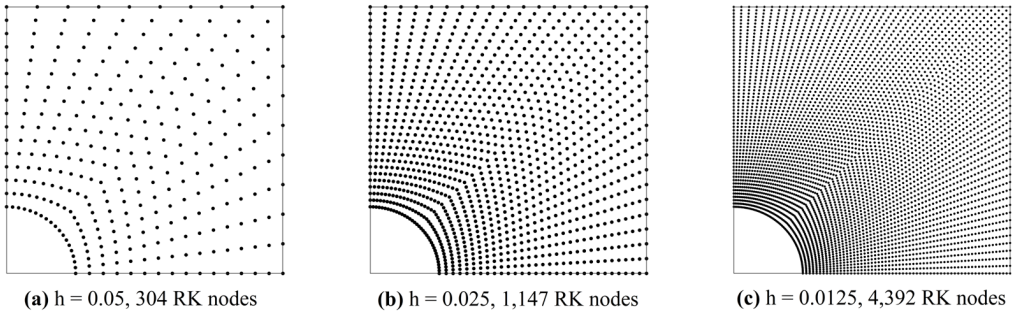

**(a)** h = 0.05, 304 RK nodes      **(b)** h = 0.025, 1,147 RK nodes      **(c)** h = 0.0125, 4,392 RK nodes

**Figure 10: Background RK discretization: (a) 304 nodes with $h = 0.05$ (b) 1,147 nodes with $h = 0.025$ and (c) 4,392 nodes with $h = 0.0125$ (h: averaged nodal spacing)**

$$\|\boldsymbol{u} - \boldsymbol{u}^h\|_0 = \sqrt{\frac{\int_\Omega \left(\boldsymbol{u}^{\text{exact}}(\boldsymbol{x}) - \boldsymbol{u}^h(\boldsymbol{x})\right)^2 d\Omega}{\int_\Omega \left(\boldsymbol{u}^{\text{exact}}(\boldsymbol{x})\right)^2 d\Omega}}. \tag{38}$$

The convergence curves for varying background RK nodal spacing ($h$) with different fixed hidden layer width ($n_{NR}$) are demonstrated in Figure 11(a), and Figure 11(b) presents the convergence curve for varying hidden layer widths for each background discretization. The NN enrichment is introduced to all discrete nodes, except for the case represented by the red curve in Figure 11(a), where only background RK approximation is considered, i.e., $u^h(\boldsymbol{x}) = u^{RK}(\boldsymbol{x})$. As indicated in Figure 11(a), the background RK solutions achieve the optimal convergence rate of 2, consistent with the analytical convergence in Eq. (30). Additionally, with a fixed background RK approximation at each domain discretization level, the NN-PU approximation yields convergence rate of 1 with respect to $\hat{n}^{1/2}$, where $\hat{n}$ is the number of neurons, again verify the error analysis result in Eq. (30).

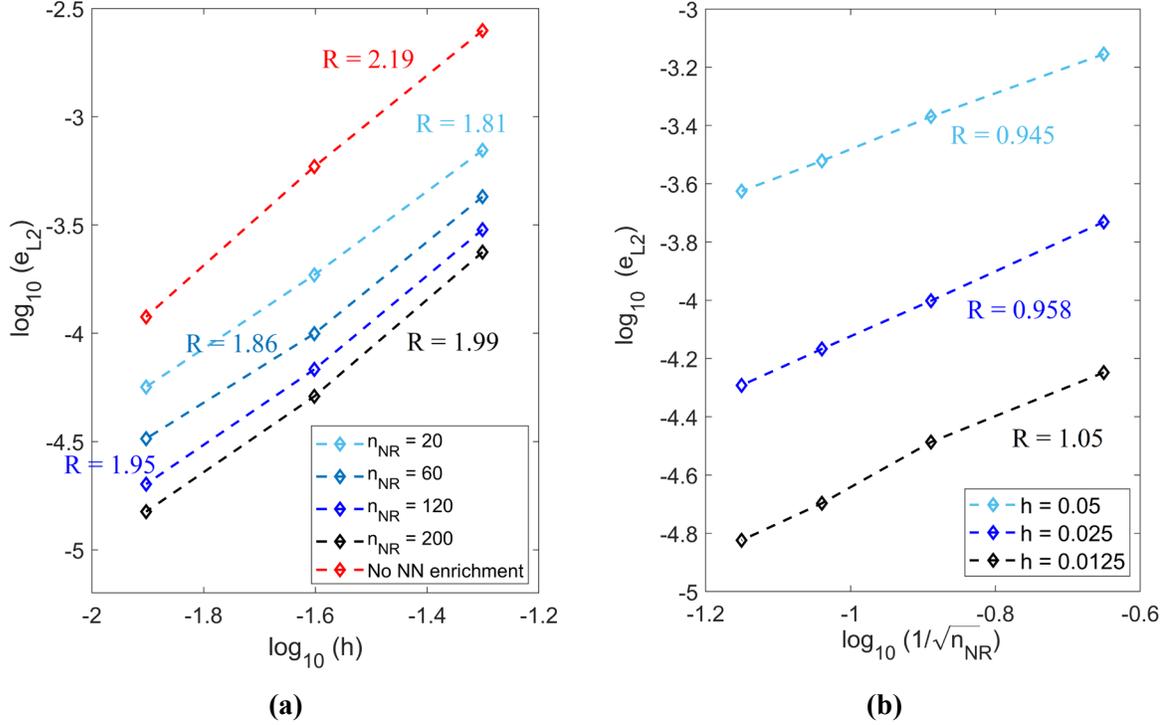

**(a)**                                    **(b)**

**Figure 11: Normalized L$_2$ convergence rates: (a) for varying background $h$ with different fixed widths of hidden layer and (b) for varying $n_{NR}$ with different level of background discretization (R: average rate of convergence)**

Next, the adaptive NN enrichment is investigated. The coarse domain discretization with $h = 0.05$, as shown in Figure 10(a) is utilized. Four cases are examined, using an energy error density ((34)-(37)) thresholds $0.1\rho_{max}^*$, $0.2\rho_{max}^*$, $0.5\rho_{max}^*$ and $0.9\rho_{max}^*$, as shown in Figure 12. We also consider enriching all nodes and adaptively selecting NN enrichment regions during the loss function minimization, as shown in Figure 13. The results of total potential energy (loss function) minimization for all cases are demonstrated in Figure 14, where the reference energy is calculated according to the analytical solutions. Enriching all nodes within the domain results in the potential energy converging closest to the analytical reference energy level among all cases. Alternatively, adaptive NN enrichment during optimization achieves comparable accuracy to that achieved by enriching all nodes.

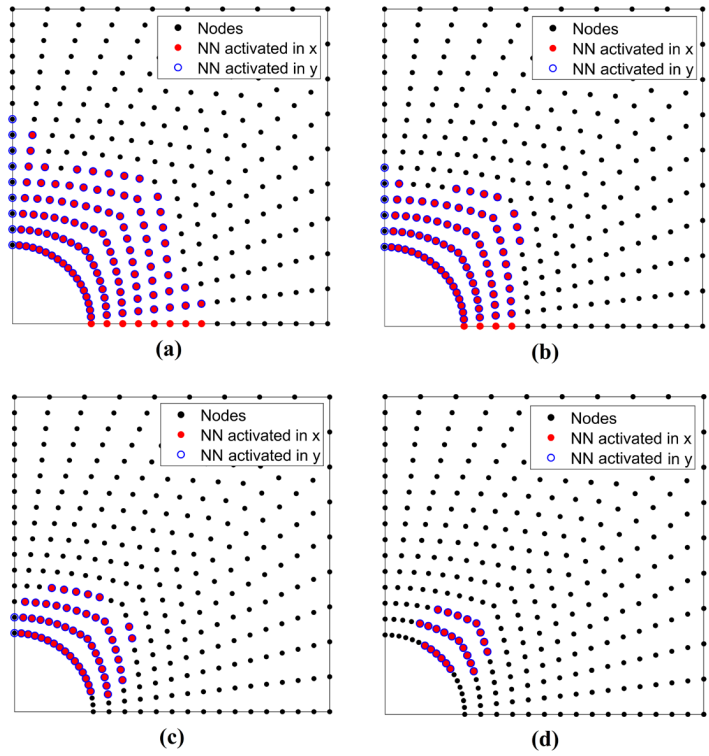

**Figure 12: NN enrichment region selections: (a) 113 and 114 nodes are enriched in x- and y- directions based on $0.1\rho_{max}^*$ (b) 80 and 82 nodes are enriched in x- and y-directions based on $0.2\rho_{max}^*$ (c) 50 and 52 nodes are enriched in x- and y-directions based on $0.5\rho_{max}^*$ and (d) 24 nodes are enriched in both directions based on $0.9\rho_{max}^*$**

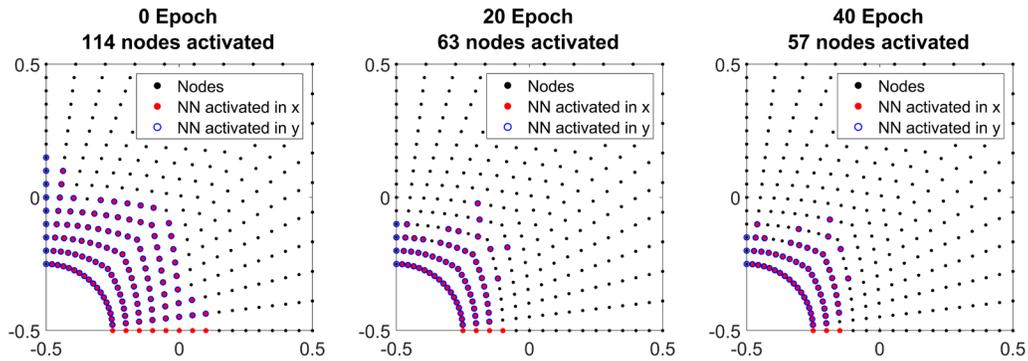

**Figure 13: Adaptive NN enrichment during loss function minimization**

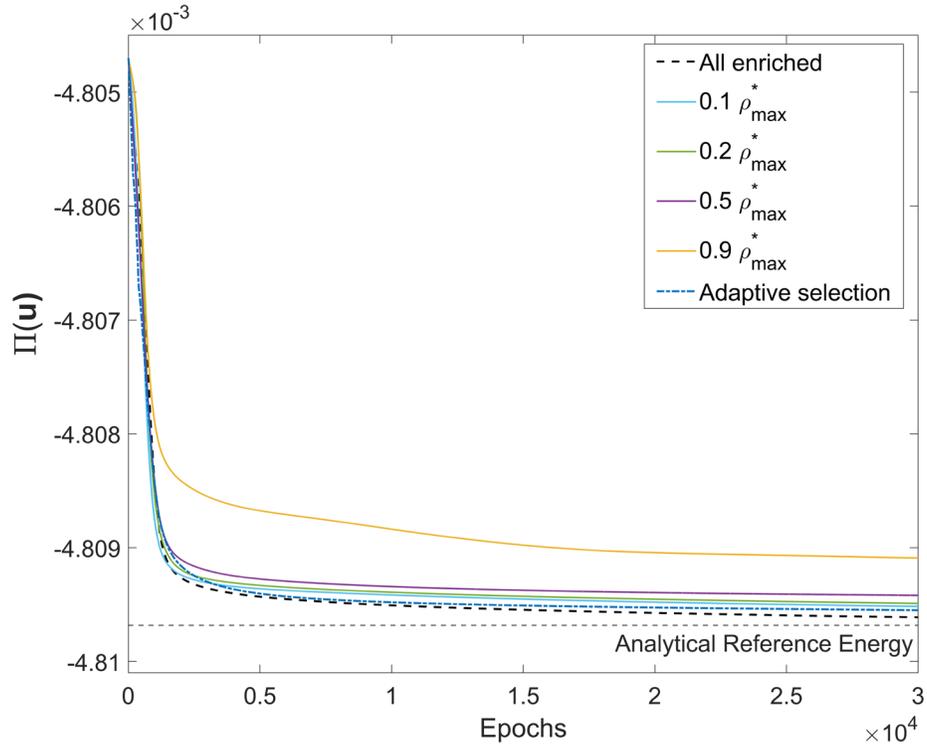

**Figure 14: Loss function minimization for NN enrichment regions using different energy error density thresholds and adaptive enrichment**

### 5.1.2 Parent Problem 2: L-shaped panel

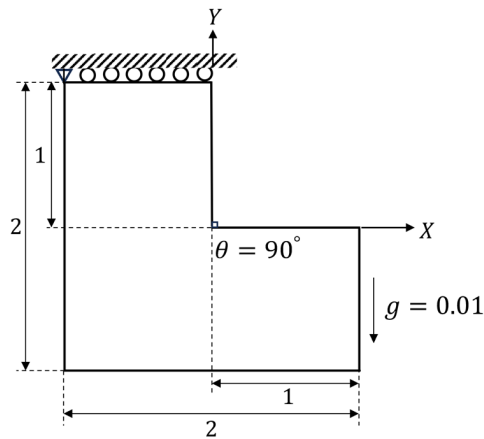

**Figure 15: L-shaped panel problem setting**

Figure 15 presents an L-shaped panel "parent" problem with a surface angle $\theta = 90°$. As a reference of the NN-PU background discretization, finite element discretization with various $h$- and $p$- refinements shown in Figure 16 are considered, and the FEM stress solutions shown in Figure 17 revealing that only the Q8 FEM solution with 578,291 nodes precisely captures the stress singularity around the concave corner.

The NN-PU background discretization is based on the nodal locations of the coarsest FE discretion with 341 nodes in Figure 16(a). Three parametric studies are performed to evaluate the effects of NN Basis sub-block's architecture on its convergence behaviors in minimizing the total potential energy functional (loss function). The initial NN-enriched nodeset $\mathcal{S}^{\zeta}$ obtained based on the algorithms in Eqs. (34)-(37) with a threshold of $0.9\rho_{max}^*$ contains 13 nodes around the concave corner, as shown in Figure 18, and they are used as the enrichment nodes for all parametric studies. A constant learning rate of $10^{-4}$ is chosen for optimization of NN enrichment basis functions.

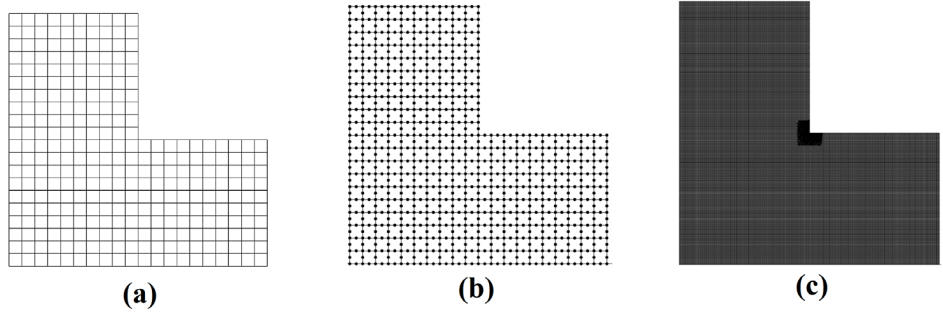

**Figure 16: Domain discretization used for obtaining FEM solutions: (a) Q4 FEM with 341 nodes (b) Q8 FEM with 981 nodes (c) Q8 FEM with 578,291 nodes**

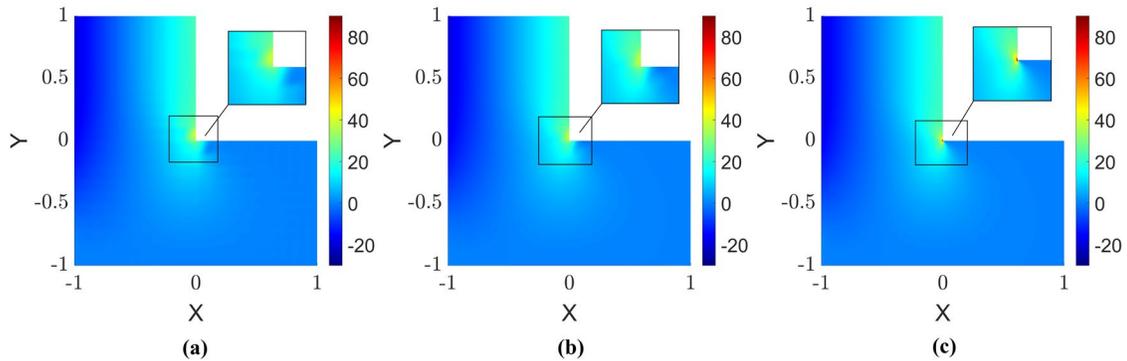

**Figure 17: Comparison of the FEM $\sigma_{22}$ results: (a) Q4 FEM with 341 nodes (b) Q8 FEM with 981 nodes (c) Q8 FEM with 578,291 nodes**

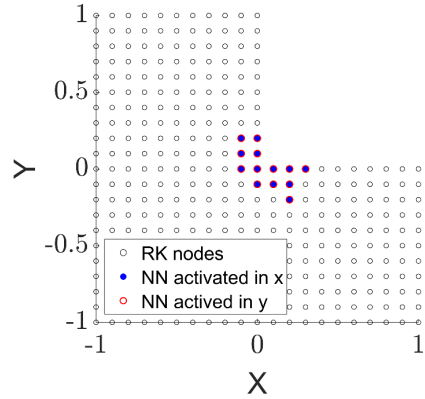

**Figure 18: Nodes activated for NN enrichment**

(1) <u>The effects of the number of NN-enrichment bases $n^\zeta$</u>: The number of NN-enrichment bases is determined by the number of neurons in the last hidden layer of NN Basis sub-block. A four-layered deep neural network is adopted for this study, and the first 3 hidden layers contain 20 neurons in each hidden layer.

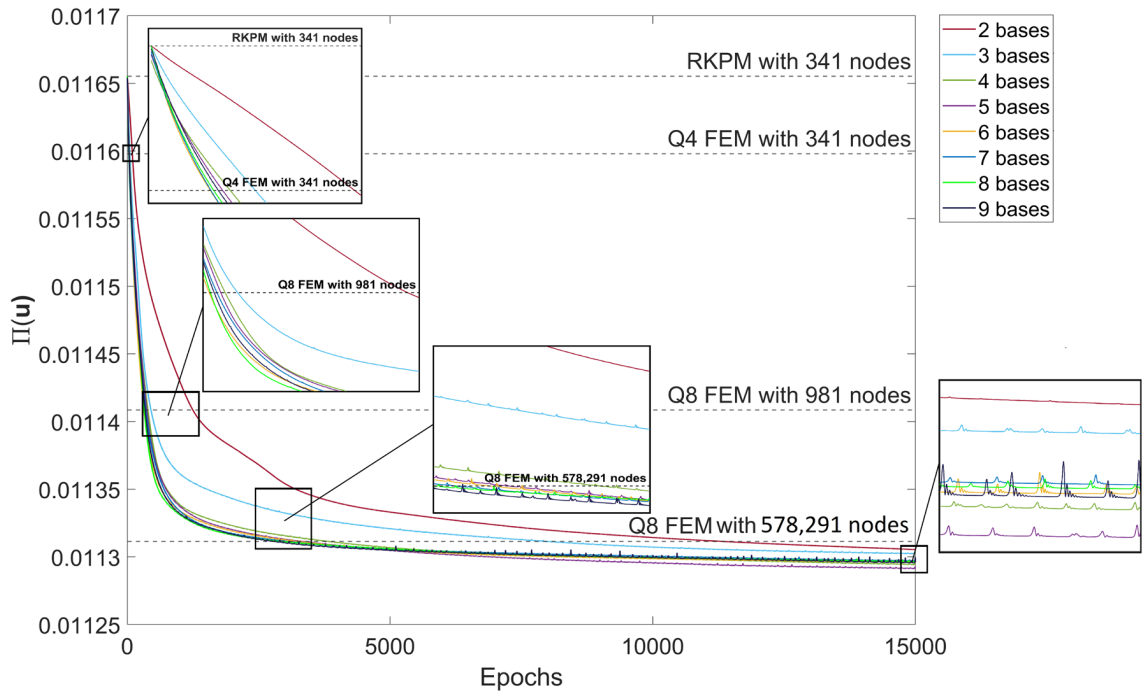

**Figure 19: Loss function (potential energy) minimization for NN architectures with different numbers of NN-enrichment bases and the comparison with FEM potential energies**

| Number of Bases | PE less than Q4 FEM with 341 nodes (PE = 1.1598E-02, CPU = 1.4s) | | PE less than Q8 FEM with 981 nodes (PE = 1.1409E-02, CPU = 2.1s) | | PE less than Q8 FEM with 578,291 nodes (PE = 1.1311E-02, CPU = 40s) | | Converged PE |
|---|---|---|---|---|---|---|---|
| | Epoch | CPU (s) | Epoch | CPU (s) | Epoch | CPU (s) | |
| 2 | 97 | 1.6 | 1,248 | 20.0 | 11,342 | 181.5 | 1.130139E-02 |
| 3 | 51 | 0.83 | 451 | 7.4 | 7,989 | 130.7 | 1.129684E-02 |
| 4 | 39 | 0.65 | 393 | 6.6 | 3,827 | 64.0 | 1.128858E-02 |
| 5 | 35 | 0.60 | 372 | 6.3 | 3,325 | 56.5 | 1.128733E-02 |
| 6 | 23 | 0.40 | 300 | 5.2 | 4,066 | 70.9 | 1.129156E-02 |
| 7 | 30 | 0.53 | 318 | 5.7 | 3,222 | 57.3 | 1.129242E-02 |
| 8 | 31 | 0.56 | 323 | 5.9 | 2,908 | 52.8 | 1.129219E-02 |
| 9 | 34 | 0.63 | 345 | 6.4 | 3,018 | 55.8 | 1.129262E-02 |

**Table 1: Summary of the potential energies (PE) and CPUs for NN architectures with different numbers of NN-enrichment bases**

It is expected that the numerical solution with lower potential energy corresponds to a higher solution accuracy. As shown in Figure 19, FE model with finer discretization and higher order of approximation generates lower potential energy. Figure 19 depicts the optimization results of loss functions (potential energies) for NN architectures featuring 2 to 9 NN-enrichment bases. Zoom-in plots display where the loss functions intersect with reference potential energy levels and the final converged potential energy levels. Table 1 summarizes the required epochs and CPUs for different NN architectures to yield potential energy right below those obtained from the FEM solutions. Increasing the number of learned NN-enrichment bases from 2 to 5 results in a noteworthy reduction in the required epochs and CPUs to achieve specific potential energy level of the much-refined FEM solution, as observed in Table 1. However, further increasing the number of NN bases beyond 5 does not lead to further reduction in the converged potential energy. Therefore, $n^\zeta$ is selected to be 5 for the remaining parametric studies in this example.

(2) The effects of number of hidden layers: In this study, each hidden layer, except for the last one, contains 20 neurons, and 5 neurons are used in the last hidden layer, which corresponds to 5 NN bases.

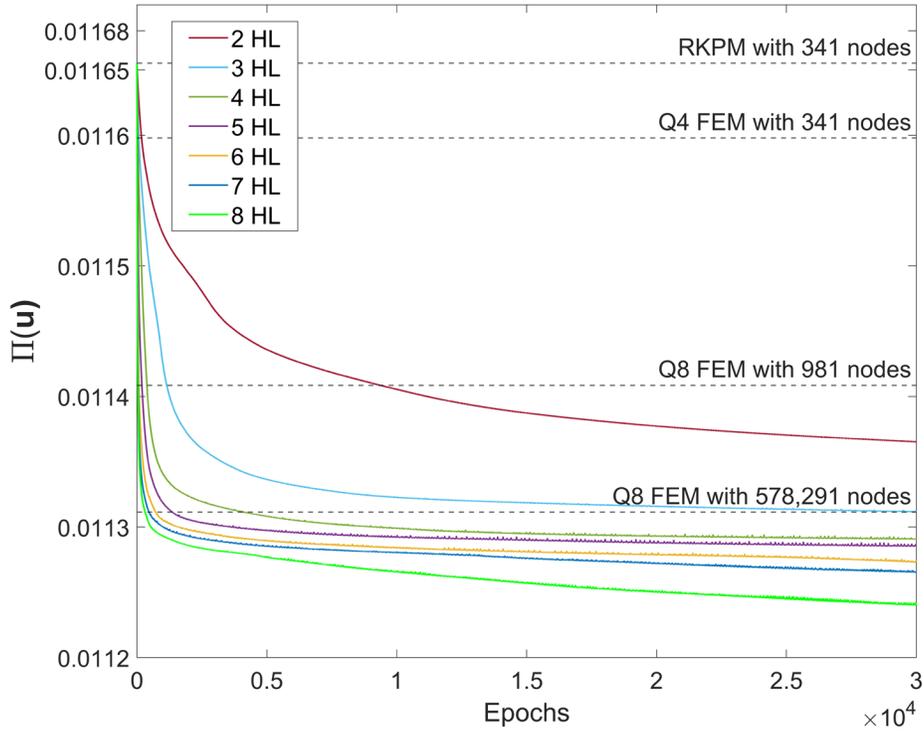

**Figure 20: Loss function minimization for NN architectures containing different number of hidden layers and the comparison with FEM potential energies**

| Number of hidden layers | PE less than Q4 FEM with 341 nodes (PE = 1.1598E-02, CPU = 1.4s) | | PE less than Q8 FEM with 981 nodes (PE = 1.1409E-02, CPU = 2.1s) | | PE less than Q8 FEM with 578,291 nodes (PE = 1.1311E-02, CPU = 40s) | | Converged PE |
|---|---|---|---|---|---|---|---|
| | Epochs | CPU (s) | Epochs | CPU (s) | Epochs | CPU (s) | |
| 2 | 178 | 1.8 | 9,360 | 93.6 | > 30,000 | -- | 1.135332E-02 |
| 3 | 73 | 0.88 | 1,146 | 17.2 | > 30,000 | -- | 1.131179E-02 |
| 4 | 35 | 0.60 | 372 | 6.3 | 3,325 | 56.5 | 1.128733E-02 |
| 5 | 21 | 0.38 | 190 | 3.4 | 1,356 | 24.4 | 1.128538E-02 |
| 6 | 11 | 0.21 | 78 | 1.5 | 716 | 13.6 | 1.127323E-02 |
| 7 | 7 | 0.14 | 48 | 0.96 | 463 | 9.3 | 1.126534E-02 |
| 8 | 4 | 0.084 | 37 | 0.78 | 313 | 6.3 | 1.125375E-02 |

**Table 2: Summary of the potential energies (PE) and CPUs for NN architectures containing different number of hidden layers**

The loss functions during the minimization process for NN architectures with 2 to 8 hidden layers are plotted in Figure 20. Table 2 summarizes the required numbers of epochs and CPUs for the loss functions of different architectures to become lower than specific reference FE potential energy levels, together with the converged potential energies. For NN Basis sub-Block comprising of 2 or 3 hidden layers, the optimization process is computationally demanding. However, it's evident that deeper NN architectures lead to significant reduction in the potential energy, as well as reduced epochs and total CPU. Employing an eight-layered NN architecture with 5 bases achieves the lowest reference potential energy level of the finest FEM solution in just 6.3 seconds, far outperforming the 40-second CPU time in FEM.

(3) <u>The effects of the number of neurons within each hidden layer</u>: The NN Basis sub-Block is constructed using 8 hidden layers and 5 bases, and the first 7 hidden layers contain an identical number of neurons.

Figure 21 demonstrates the potential energies for NN architectures featuring 10 to 80 neurons per hidden layer. Table 3 documents the epochs and CPUs required for various NN architectures to generate potential energy below the reference FEM potential energies and their respective converged energy levels. The results clearly indicate that increasing the width of the hidden layer from 10 neurons to 40 neurons yield a notable acceleration in the optimization speed of the neural network, as the NN-enrichment bases contain more trainable parameters to capture the local features. However, as the number of neurons per hidden layer further increases, the improvement becomes marginal. Overall, when employing 8 hidden layers and 5 bases for constructing the NN-enrichment, the offline learning is shown to be more efficient than the Q8 FEM solution with 579,012 nodes (FEM CPU = 40 Sec), especially when the hidden layer has a width of 20 neurons or more.

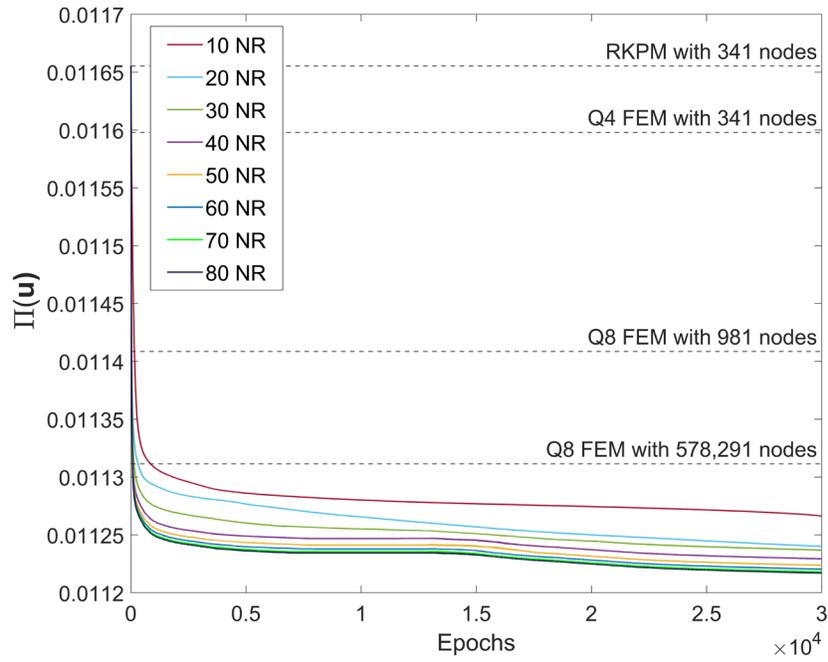

**Figure 21: Loss function minimization for NN architectures containing different number of neurons per hidden layer and the comparison with FEM potential energies**

| No. of Neurons | PE less than Q4 FEM with 341 nodes (PE = 1.1598E-02, CPU = 1.4s) | | PE less than Q8 FEM with 981 nodes (PE = 1.1409E-02, CPU = 2.1s) | | PE less than Q8 FEM with 578,291 nodes (PE = 1.1311E-02, CPU = 40s) | | Converged PE |
|---|---|---|---|---|---|---|---|
| | Epochs | CPU (s) | Epochs | CPU (s) | Epochs | CPU (s) | |
| 10 | 22 | 0.46 | 155 | 3.26 | 860 | `18.1 | 1.126635E-02 |
| 20 | 4 | 0.08 | 37 | 0.78 | 313 | 6.3 | 1.125375E-02 |
| 30 | 5 | 0.11 | 32 | 0.70 | 140 | 3.1 | 1.123675E-02 |
| 40 | 5 | 0.12 | 29 | 0.67 | 86 | 2.0 | 1.122940E-02 |
| 50 | 5 | 0.13 | 29 | 0.73 | 77 | 1.9 | 1.122378E-02 |
| 60 | 6 | 0.16 | 31 | 0.81 | 72 | 1.9 | 1.122042E-02 |
| 70 | 7 | 0.20 | 31 | 0.90 | 70 | 2.0 | 1.121817E-02 |
| 80 | 6 | 0.19 | 32 | 0.99 | 71 | 2.2 | 1.121705E-02 |

**Table 3: Summary of the potential energies (PE) and CPUs for NN architectures containing different number of neurons per hidden layer**

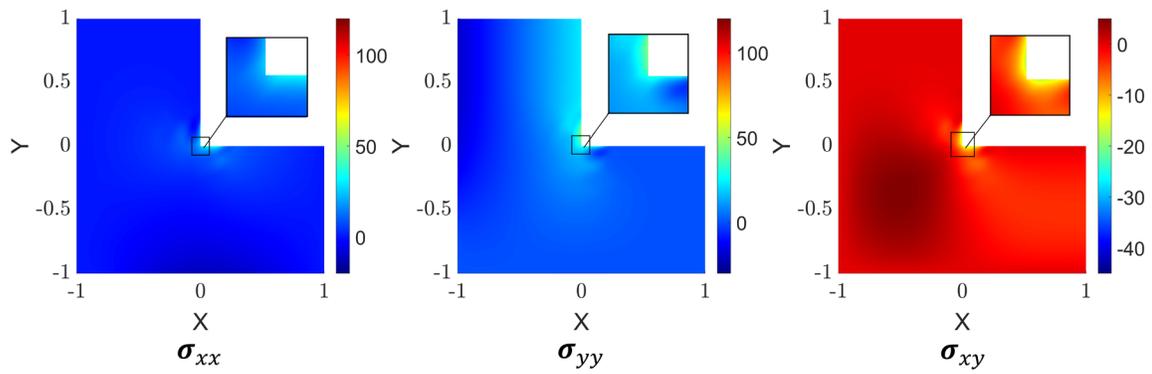

**(a) RKPM with 341 nodes**

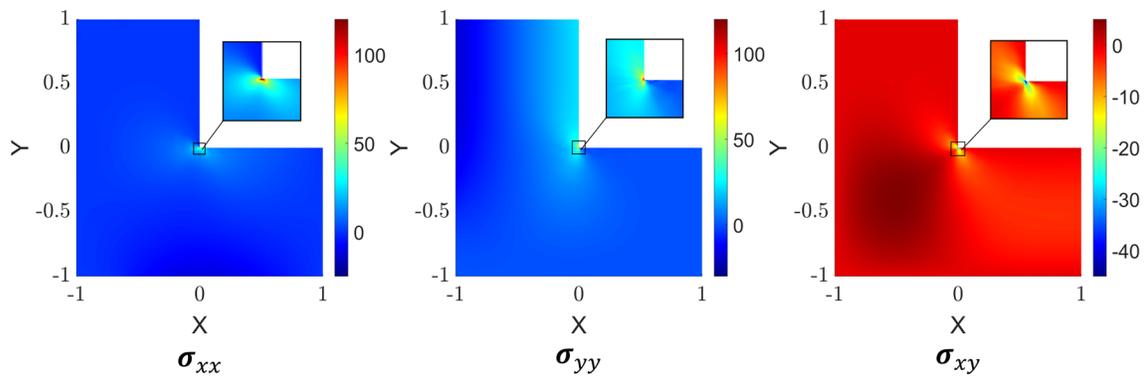

**(b) NN-PU with 341 RK nodes and 2,815 NN unknowns**

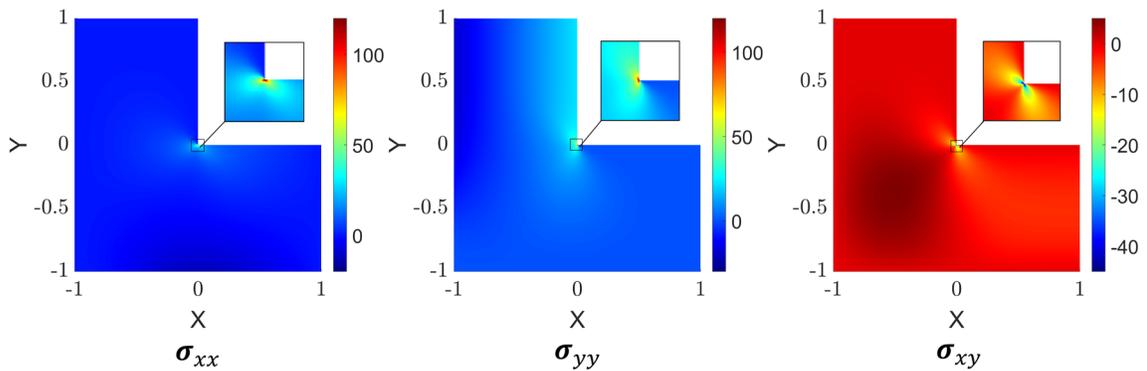

**(c) Q8 FEM with 578,291 nodes**

**Figure 22: Comparison of the stress solutions: (a) RKPM solutions with 341 nodes**
**(b) NN-PU solutions with 341 RK nodes and 2,815 NN unknowns and (c) Q8**
**FEM solutions with 578,291 nodes**

Figure 22 demonstrates a comparison between the stress distributions obtained using the RKPM with 341 nodes, NN-PU approximaotion with 341 nodes and 2,815 NN hyperparameters, and Q8 FEM with 578,291 nodes (over a million degrees of freedom). The NN Basis sub-block comprises 8 hidden layers with 20 neurons each, and 5 bases, totaling 2,685 trainable parameters. The NN enrichment is applied to 13 nodes, each with 2 degrees of freedom, yielding 130 NN enrichment coefficients. Consequently, the total number of NN unknowns is 2,815. As shown in Figure 22, the NN-PU approximation is able to obtain indistinguishable solutions compared to the results obtained using the Q8 FEM with a much finer discretization with just 0.3% DOFs and 16.6% CPU of the FEM computation.

## 5.2 Online analysis with transfer learning

This subsection demonstrates how the NN basis functions constructed for local features in the "parent" problems introduced in Sec. 5.1 can be used as the NN enrichments for problems with the same or similar local features via a transfer learning strategy as introduced in Sec. 4.1.2.

### 5.2.1 Circular-shaped panel with a 90-degree corner: transfer learned from the L-shaped parent problem

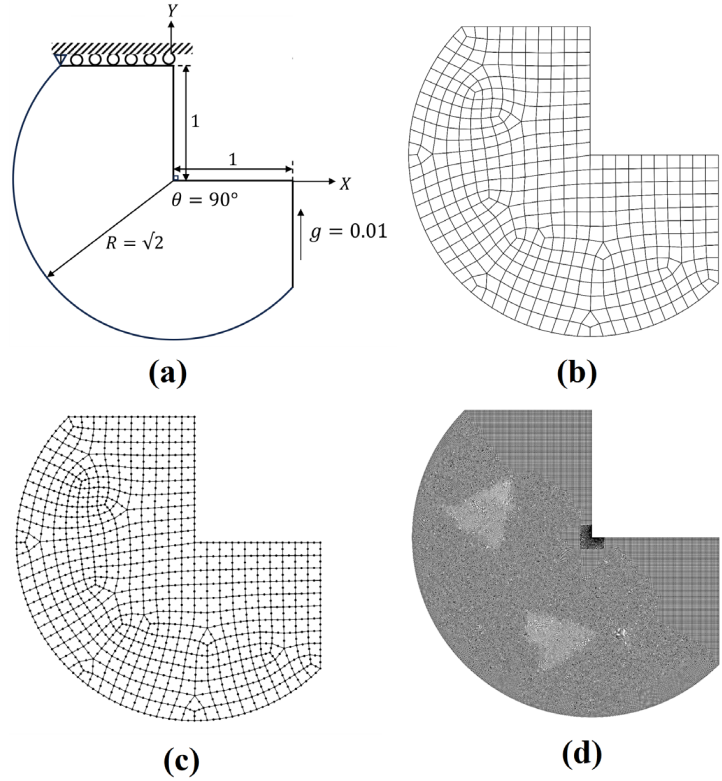

**(a)**  **(b)**

**(c)**  **(d)**

**Figure 23: (a) L-shape panel problem setting with a circular domain and 3 levels of FE discretization: (b) Q4 FEM with 480 nodes (c) Q8 FEM with 1,396 nodes (d) Q8 FEM with 776,438 nodes**

Figure 23 presents a variation of the "parent" L-shaped panel discussed in Sec. 5.1.2 to a circular-shaped panel with a 90-degree corner. Despite of these domain geometry changes, the problem maintains its local 90-degree corner same as that in the L-shaped "parent" problem. Utilizing a "parent" NN-enrichment basis set obtained from the L-shaped problem, transfer learning is applied to the online solution stage of this circular domain problem. We incorporate a Parametric sub-block $\mathcal{N}^P$ comprising 2 hidden layers with 5 neurons in the first hidden layer, 2 neurons in the second hidden layer, totaling of 27 trainable parameters, to facilitate the utilization of the pre-trained NN enrichment basis sets during online computations. The same set of NN basis functions in the L-shaped problem are employed without additional optimization in their spatial distributions, and only the weight sets of the Parametric sub-block ($\mathbf{W}^P$) and NN basis function coefficients ($\mathbf{W}^C$) are iterated in the loss function minimization.

The chosen NN basis sub-block architecture from the "parent" L-shaped problem is composed of 8 hidden layers, with the first 7 containing 20 neurons each, and it considers 5 NN-enrichment bases as studied in Sec. 5.1.2. During the online solution stage, the NN enriched nodeset is updated every 20 epochs with a threshold of $0.1\rho_{max}^*$ until the number of enriched nodes remains unchanged. Figure 24 demonstrates the adaptive enrichment process of the NN-enriched nodes during the optimization, where the neural network efficiently identifies the areas where enrichment is most needed within just 40 epochs, consolidating the number of enriched nodes from 180 to 29.

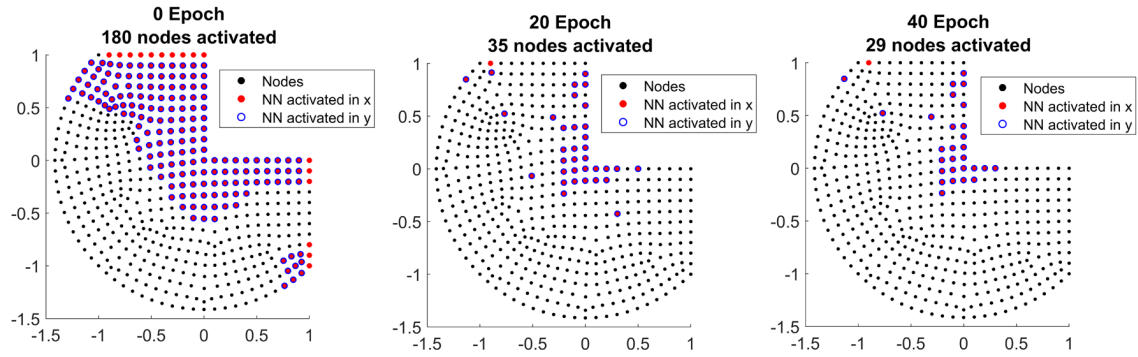

**Figure 24: Adaptive selection of the nodes for NN enrichment**

Figure 25 demonstrates the NN-PU loss function minimization in the online solution stage in comparison with various reference potential energy levels computed using FEM solutions with varying spatial discretization levels shown in Figure 23. Table 4 provides a detailed comparison of the number of unknowns, the discrete total potential energy levels and computation time required for NN-PU and FEM with different domain discretization. In addition, it also reports the number of required epochs and CPU for the NN-PU online calculation with the learned "parent" NN-enrichment basis set to reach each FEM potential energy level. When compared to the high-order highly refined FEM model (with 1,551,873 DOFs), the NN-PU approximation achieves a similar level of solution accuracy with just 317 NN variables (1,254 DOFs total) as shown in Figure 26, while requiring only 20.4% of computation time in the online computation as shown in Table 4. The final converged potential energy for the online solution is 1.6518E-02, 3.74 % lower (more accurate) than that from Q8 FEM with 776,438 nodes (1,551,873 DOFs).

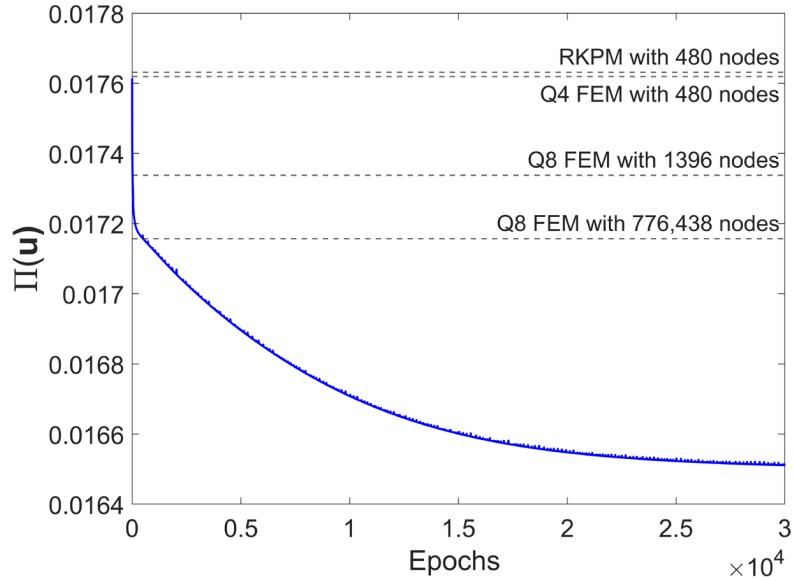

**Figure 25: Loss function minimization of NN-PU and the comparison with FEM potential energies (PEs)**

| | Q4 FEM with 480 nodes (PE = 1.762E-02, CPU = 1.9s) | Q8 FEM with 1,396 nodes (PE = 1.734E-02, CPU = 2.9s) | Q8 FEM with 776,438 nodes (PE = 1.716E-02, CPU = 56.3s) |
|---|---|---|---|
| NN-PU with potential energy (PE) lower than FEM PE | | | |
| Number of Epochs NN-PU with 480 nodes used to reach each FEM potential energy levels | 1 | 16 | 468 |
| CPU for NN-PU with 480 nodes to reach each FEM potential energy levels (s) | 2.2 | 2.4 | 11.5 |

**Table 4: Comparison of the required numbers of epochs and runtimes for NN-PU solution to reach the potential energy (PE) of FEM with different discretization for the circular L-shaped panel problem**

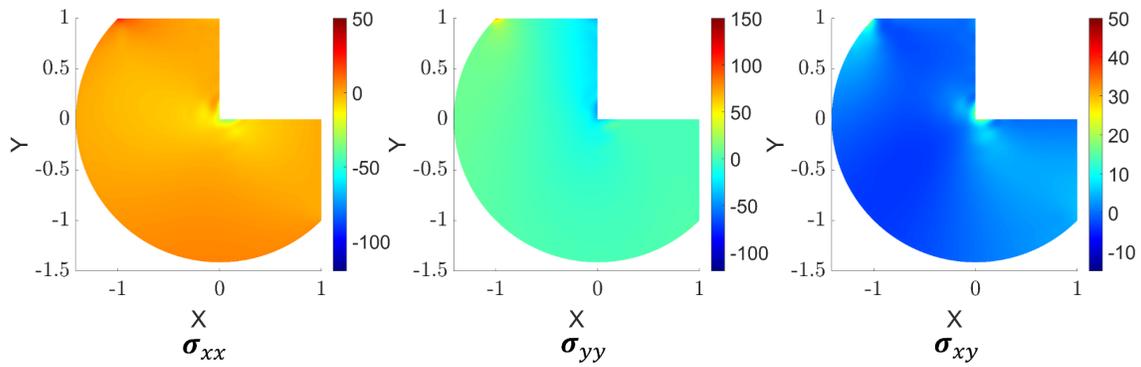

**(a) RKPM with 480 nodes**

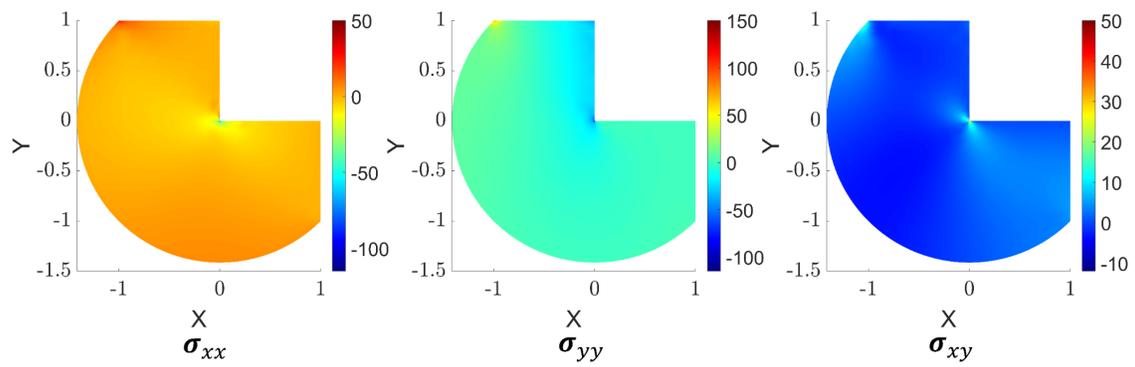

**(b) NN-PU with 480 RK nodes and 317 NN unknowns**

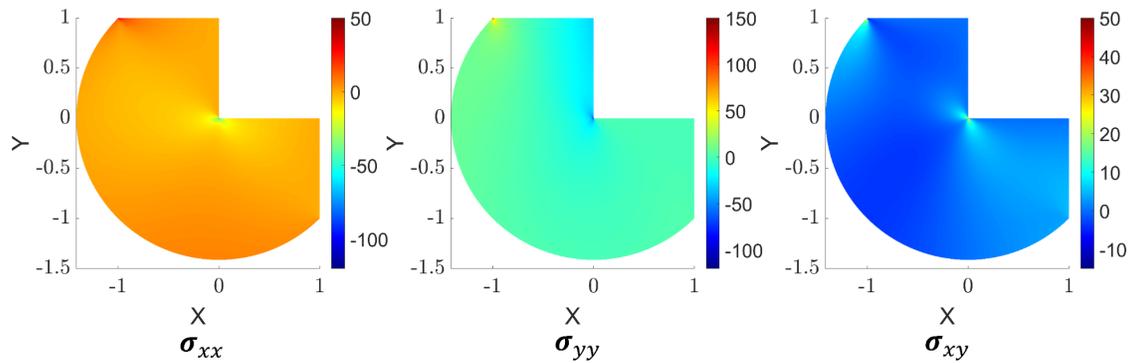

**(c) Q8 FEM with 776,438 nodes**

**Figure 26: Comparison of the stress solutions: (a) RKPM solutions with 480 nodes (b) NN-PU solutions with 480 RK nodes and 317 NN unknowns and (c) Q8 FEM solutions with 776,438 nodes**

### 5.2.2 L-shape panel with 60° corner: transfer learned from the L-shaped parent problem

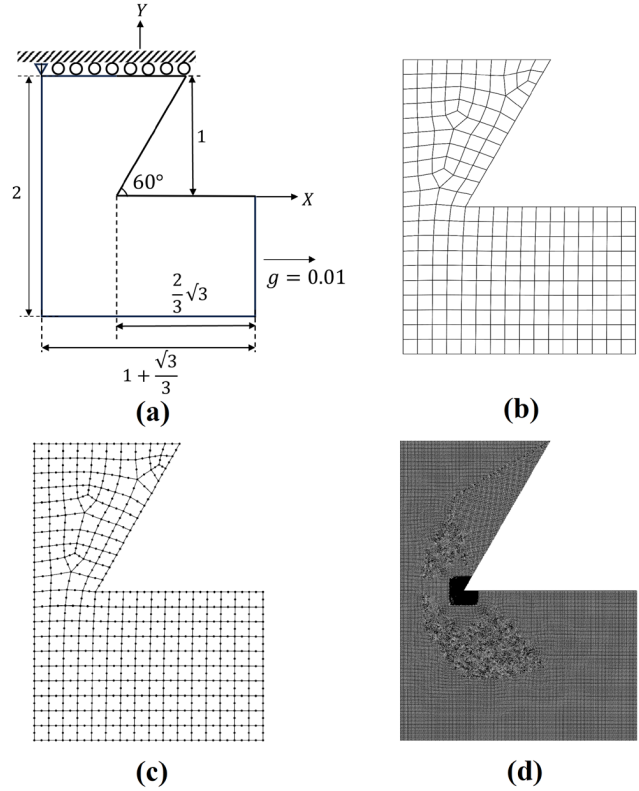

**(a)**

**(b)**

**(c)**

**(d)**

**Figure 27: (a) L-shape panel with 60° corner problem setting and 3 levels of FE discretization: (b) Q4 FEM with 273 nodes (c) Q8 FEM with 777 nodes (d) Q8 FEM with 447,602 nodes**

Here we consider a L-shaped panel problem featuring a 60° interior corner as shown in Figure 27. Different from the "parent" L-shaped panel with 90° corner, the degree of singularity around the concave corner has changed, and the essential boundary condition is now applied along the positive x-axis. Due to the change of local feature, we use a larger Parametric sub-block $\mathcal{N}^P$ comprising 3 hidden layers with 20 neurons in the first two hidden layers and 2 neurons in the last hidden layer, totaling of 522 trainable parameters. The NN Basis sub-block is transfer learned from the 90° L-shaped problem in Sec. 5.1.2 again containing 8 hidden layers, with 20 neurons in the first 7 layers, and 5 NN-enrichment bases in the last hidden layer. Like the previous problems, the selection of NN-enriched nodes is dynamically adjusted during the loss function minimization process with an activation threshold of $0.1\rho^*_{max}$. As shown in Figure 28, the initial set of 148 nodes chosen for NN-enrichment is subsequently condensed to just 50 nodes situated around the corner.

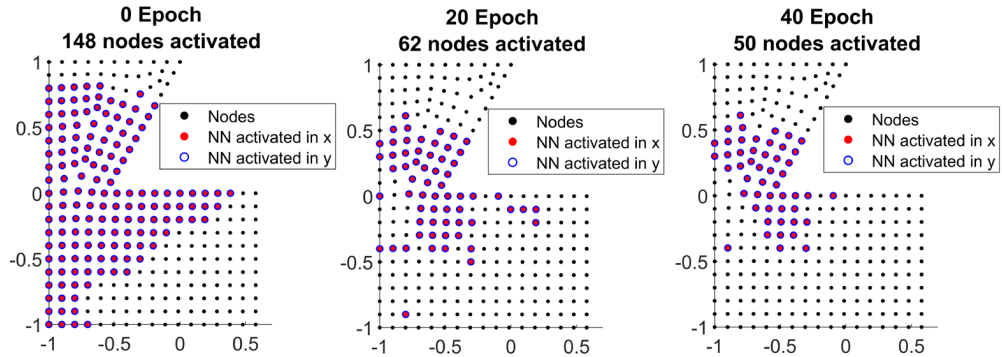

**Figure 28: Adaptive selection of the nodes for NN enrichment**

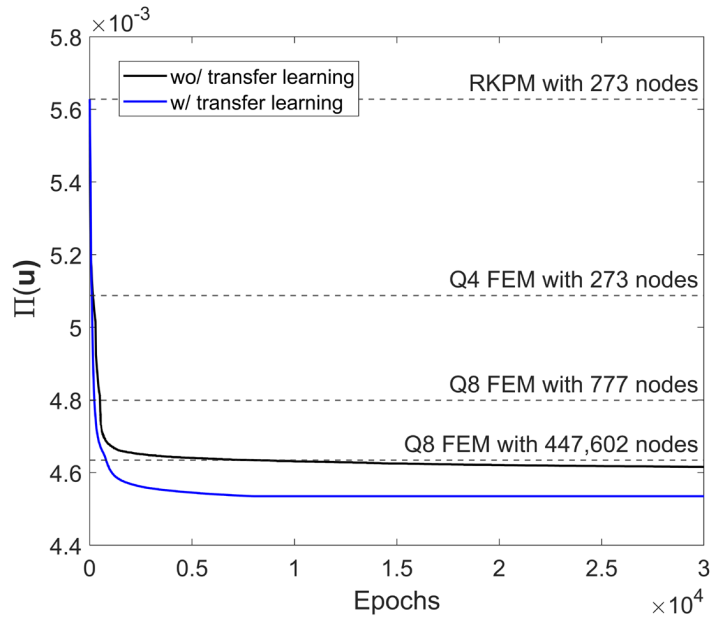

**Figure 29: Loss function minimization for NN-PU with and without transfer learning and the comparison with FEM potential energies**

Figure 29 demonstrates the loss function reduction during the optimization process with and without a transfer learning from pre-trained "parent" NN enrichment basis set. Table 5 compares the epochs and CPU required for NN-PU solution computation with and without transfer learning to yield the potential energy levels of FEM approximations with different domain discretization. The results show that utilizing transfer learned NN Basis sub-block weight set notably accelerates the optimization process compared to the case without transfer learning where random number initialization of NN Basis sub-block weights and biases were considered. Figure 30 presents the stress solutions obtained using three different approaches: RKPM with 273 nodes, NN-PU approximation with the transfer-learned NN enrichment bases, and Q8 FEM with 447,602 nodes. While RKPM without NN enrichment exhibits oscillations near the corner due to the employment of SCNI [54] without stabilization, NN-PU effectively mitigates these oscillations without any stabilization. Compared to the high-order much refined FEM solution with nearly a million degrees of freedom, NN-PU with the transfer-learned NN-enrichment basis achieves comparable accuracy with less than 0.18% of the degrees of freedom and reduced computation time. The final converged potential energy for the NN-PU online solution is 4.5350E-03, leading to a 2.14% potential energy reduction (better accuracy) compared to the one from the Q8 FEM solution with 447,602 nodes. The results demonstrate that the proposed NN-PU with transfer learning from parent problem with different features becomes more effective than FEM when higher solution accuracy is in demand.

| NN-PU with potential energy (PE) lower than FEM PE | | Q4 FEM with 273 nodes (PE = 5.088E-03, CPU = 1.3s) | Q8 FEM with 777 nodes (PE = 4.799E-03, CPU = 1.9s) | Q8 FEM with 447,602 nodes (PE = 4.634E-03, CPU = 32.9s) |
|---|---|---|---|---|
| Number of Epochs to reach the accuracy | w/ transfer learning | 88 | 201 | 802 |
| | wo/ transfer learning | 160 | 501 | 8018 |
| CPU for NN-PU (s) | w/ transfer learning | 4.5 | 6.9 | 20.2 |
| | wo/ transfer learning | 6.0 | 13.5 | 178.9 |

**Table 5: Comparison of the required numbers of epochs and runtimes for NN-PU solution to reach the potential energy (PE) of FEM solutions with different discretization with and without weights transfer learned from the parent NN basis set for the 60° corner problem**

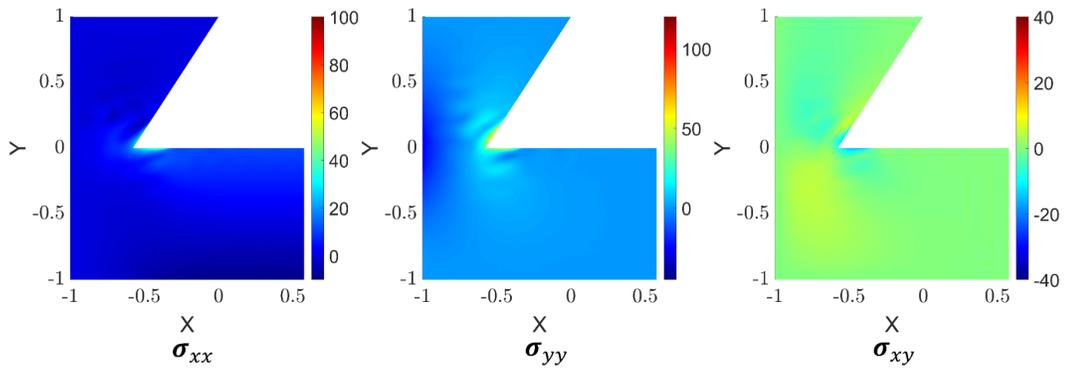

**(a) RKPM with 273 nodes**

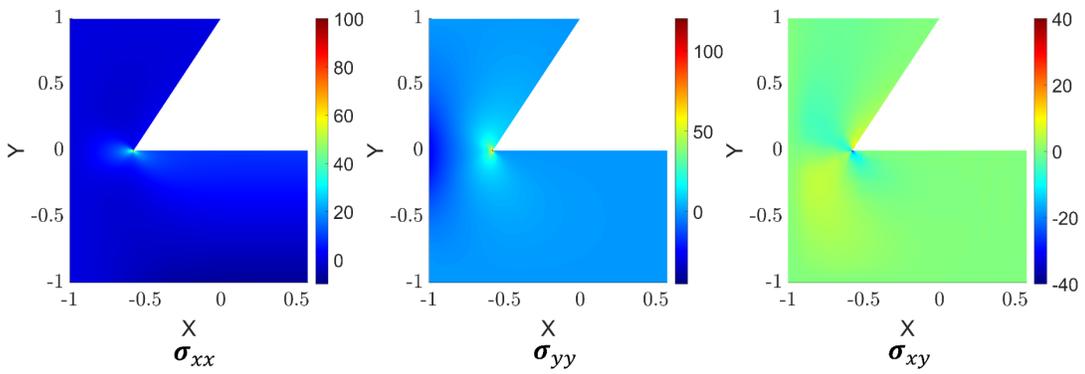

**(b) NN-PU with 273 RK nodes and 1,022 NN unknowns**

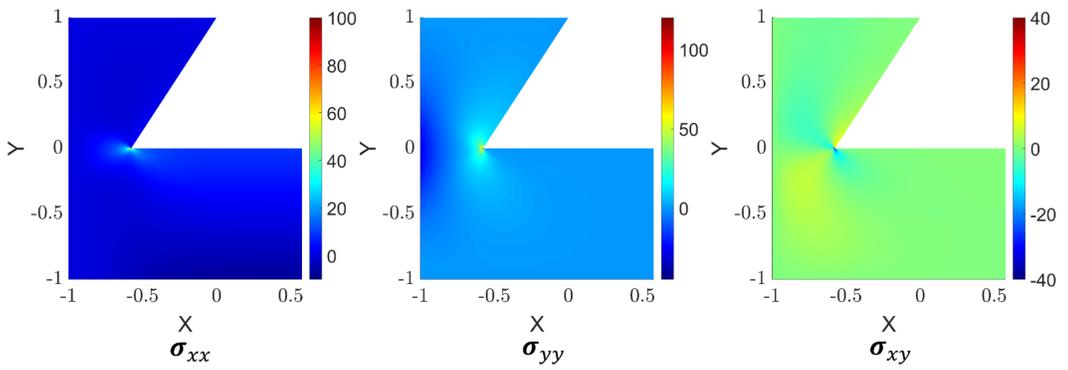

**(c) Q8 FEM with 447,602 nodes**

**Figure 30: Comparison of the stress solutions: (a) RKPM solutions with 273 nodes (b) NN-PU solutions with 273 RK nodes and 1,022 NN unknowns and (c) Q8 FEM solutions with 447,602 nodes**

## 5.3    Transfer learning with multiple local features

### 5.3.1    Plate with two circular holes: transfer learned from the plate with a hole parent problem

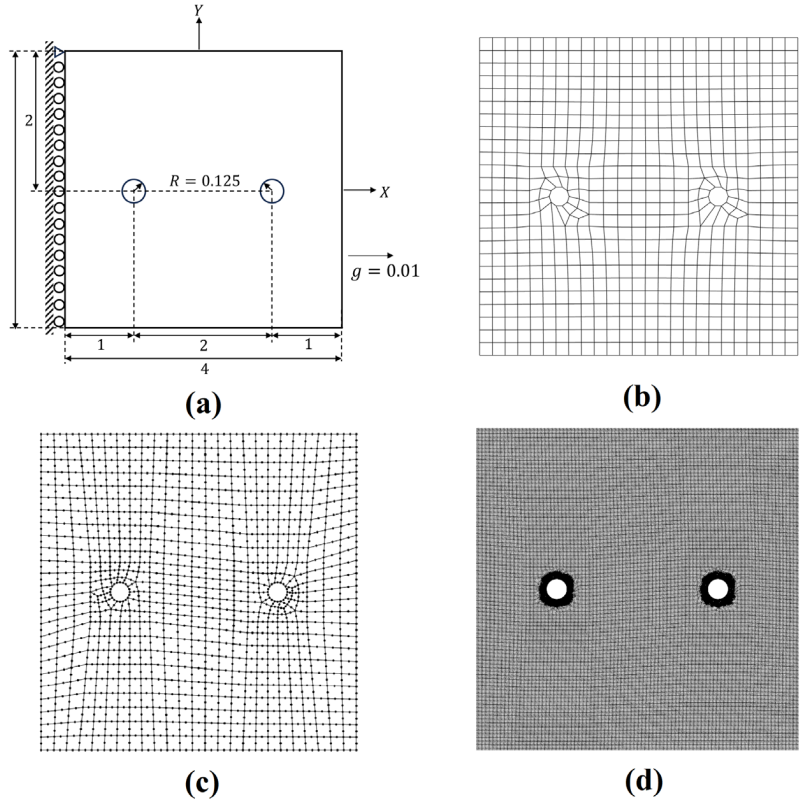

(a)                    (b)

(c)                    (d)

**Figure 31: (a) Problem settings of plate with two small circular holes and 3 levels of FE discretization: (b) Q4 FEM with 688 nodes (c) Q8 FEM with 2,015 nodes (d) Q8 FEM with 499,843 nodes**

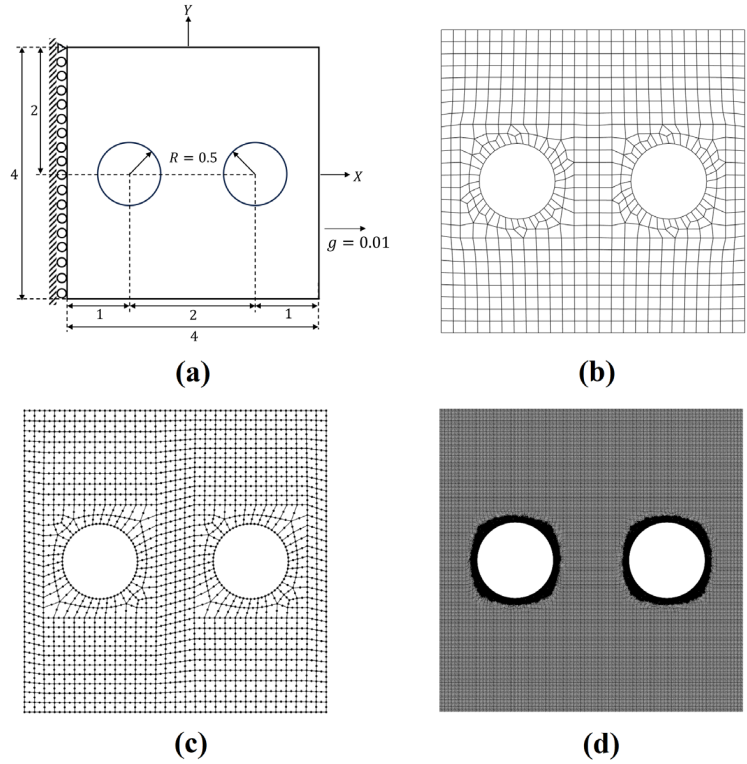

**Figure 32: (a) Problem settings of plate with two large circular holes and 3 levels of FE discretization: (b) Q4 FEM with 728 nodes (c) Q8 FEM with 2,625 nodes (d) Q8 FEM with 504,452 nodes**

In this example we consider a plate with two identical circular holes subjected to a x-directional uniaxial tension, as illustrated in Figure 31 and Figure 32, where the only difference between the two lies in the radius of the holes. Given such a problem setting, the same localized stress pattern is anticipated to be present around both circular holes. Two NN blocks are utilized for this problem for both plate geometries, and two Parametric sub-blocks are employed, denoted as $\mathcal{N}_1^P$ and $\mathcal{N}_2^P$. Both Parametric sub-blocks use the same architecture as mentioned in Sec. 5.2.1. The NN Basis sub-blocks $\mathcal{N}_1^\zeta$ and $\mathcal{N}_2^\zeta$ incorporate the same learned NN-enrichment basis functions from the plate-with-a-hole "parent" problem discussed in Sec. 5.1.1. An NN basis sub-block consisting of 4 hidden layers, with the first 3 layers comprising 50 neurons each and incorporating 5 NN-enrichment bases in the last hidden layer learned from the "parent" problem in Sec. 5.1.1 is adapted. Note that all weights and biases in $\mathcal{N}_1^\zeta$ and $\mathcal{N}_2^\zeta$ are directly transferred from those in the "parent" problem, and hence only the weight sets of Parametric sub-blocks ($\mathbb{W}^P$) and the coefficients of the NN basis functions ($\mathbb{W}^C$) are iterated during the loss function minimization.

Note that every neural network block is associated with its unique NN-enriched nodeset, referred as $\mathcal{S}_1^\zeta$ and $\mathcal{S}_2^\zeta$, respectively, which are initialized through an error estimator as follows:

$$\mathcal{S}_1^\zeta = \{I | x_{1I} < 0, \rho_I^* > 0.1\rho_{max}^*\},$$
$$\mathcal{S}_2^\zeta = \{I | x_{1I} \geq 0, \rho_I^* > 0.1\rho_{max}^*\}, \tag{39}$$

where $x_{1I}$ is the $x$ directional coordinate of a node $I$. Similar to previous problems, the nodal energy error densities are re-evaluated every 20 epochs, and both $\mathcal{S}_1^\zeta$ and $\mathcal{S}_2^\zeta$ are updated at the online computation. Figure 33 and Figure 34 demonstrate the evolutions of $\mathcal{S}_1^\zeta$ and $\mathcal{S}_2^\zeta$ for both cases, and the numbers of NN-enriched nodes are settled after three updates.

Figure 35 illustrates the process of loss function minimization for both plate geometries. When examining a plate with larger holes, as seen Figure 35(b), it becomes evident that the optimization process leads to faster convergence due to smaller interactions between the local solutions near two holes. Table 6 reports the results from NN-PU and FEM with different levels of domain discretization and Figure 36 and Figure 37 compare the stress results for plate with small holes and larger holes problems, respectively, obtained using coarse RKPM, NN-PU approximation, and high-order highly refined FEM with body-fitted meshes. The results indicate that when employing the "parent" NN-enriched enrichment basis set, NN-PU with less than 1% of the degrees of freedom can achieve similar accuracy to FEM with substantial CPU reduction in the online calculation. This efficiency is particularly evident in the case of larger holes, where the runtime is 1:7.25 compared to that of the finest FEM model. The final converged total potential energy is 5.3921E-01 for the plate with smaller holes and 4.4167E-01 for the plate with larger holes, leading to potential energy reductions of 0.08% and 0.05%, respectively, compared to the corresponding values computed using Q8 FEM with much refined meshes.

The corresponding CPUs are 34.0% and 13.8% of the corresponding FEM CPUs, respectively.

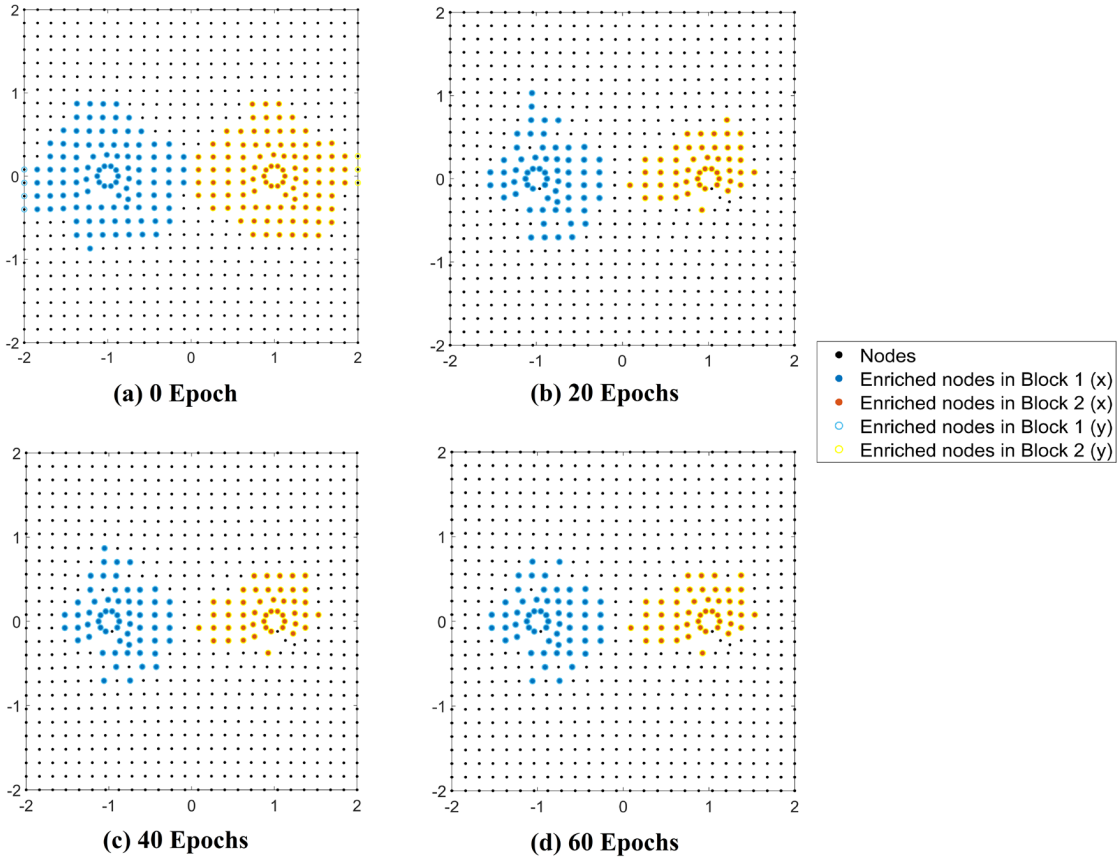

**Figure 33: Adaptive selection of the NN-enriched nodes for the problem with small circular holes: (a) total 206 and 213 enriched nodes in x- and y-directions (b) total 133 enriched nodes in both directions (c) total 109 enriched nodes in both directions and (d) total 105 enriched nodes in both directions**

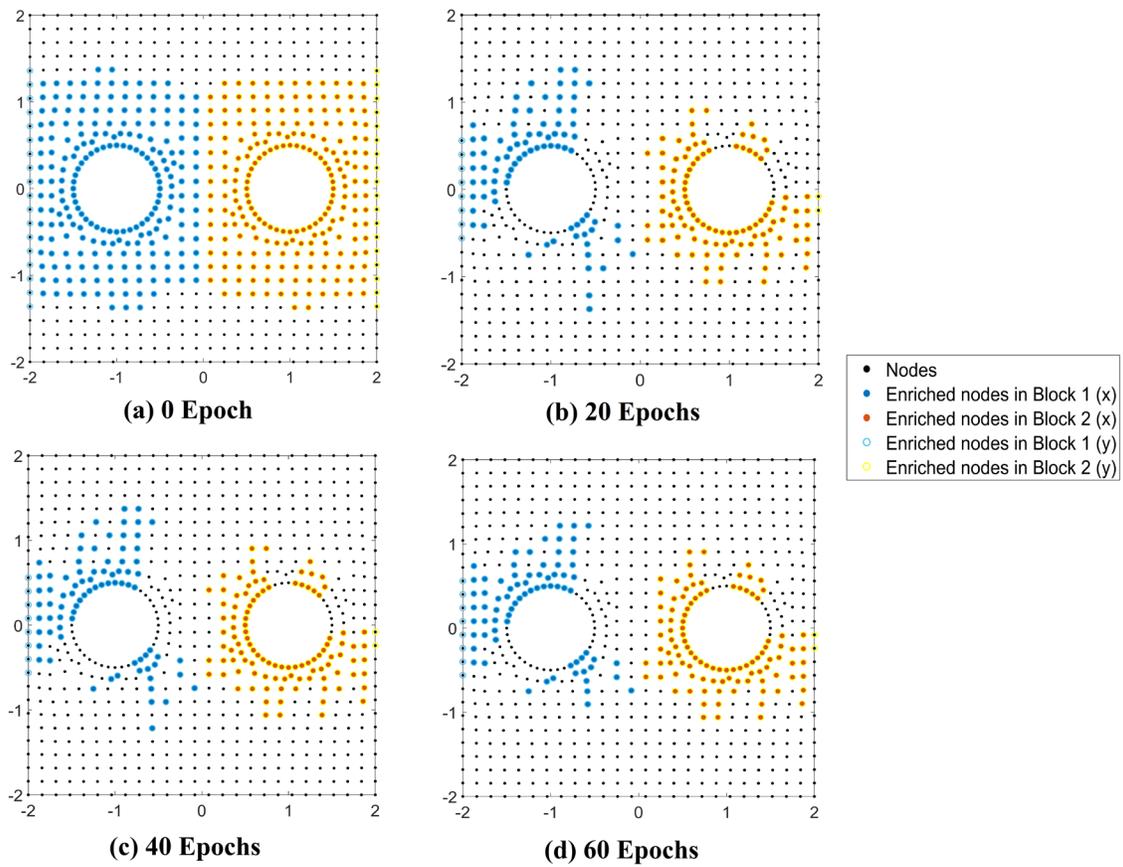

**Figure 34: Adaptive selection of the NN-enriched nodes for the problem with large circular holes: (a) total 436 and 472 enriched nodes in x- and y-directions (b) total 188 and 198 enriched nodes in x- and y-directions (c) total 187 and 197 enriched nodes in x- and y-directions and (d) total 180 and 190 enriched nodes in x- and y-directions**

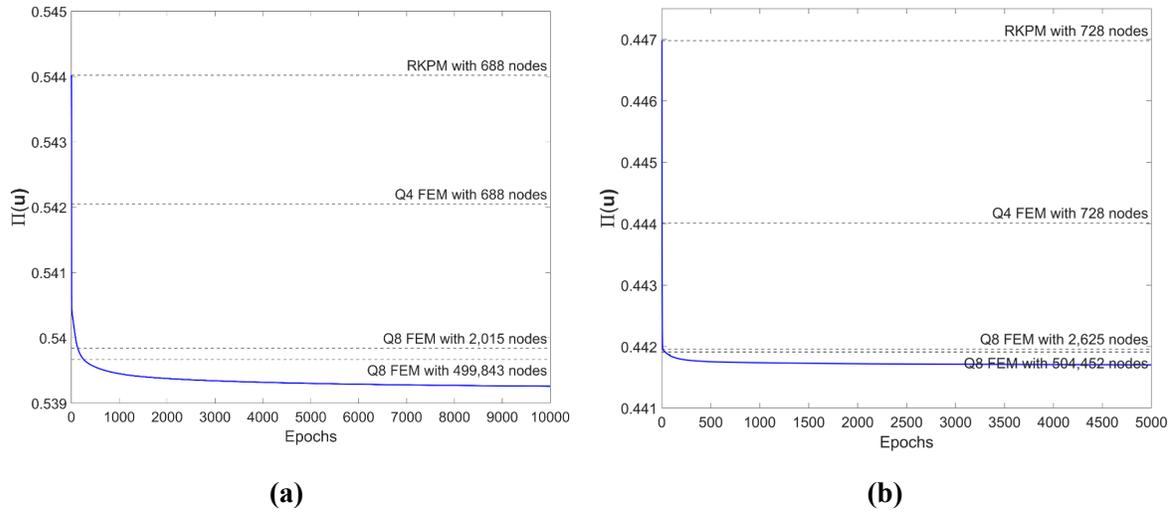

**(a)**                                    **(b)**

**Figure 35: Loss function minimization and the comparison with FEM potential energies (a): plate with two small circular holes (b) plate with two large holes**

| Two small holes with R = 0.125 | | | |
|---|---|---|---|
| NN-PU with potential energy (PE) lower than FEM PE | Q4 FEM with 688 nodes (PE = 5.421E-01, CPU=1.8s) | Q8 FEM with 2,015 nodes (PE = 5.398E-01, CPU= 2.1s) | Q8 FEM with 499,843 nodes (PE = 5.397E-01, CPU=30.3s) |
| Number of Epochs NN-PU used to reach the FEM PE levels | 3 | 135 | 257 |
| CPU for NN-PU to reach FEM PE levels (s) | 2.7 | 6.6 | 10.3 |
| Two large holes with R = 0.5 | | | |
| NN-PU with potential energy (PE) lower than FEM PE | Q4 FEM with 728 nodes (PE = 4.440E-01, CPU=1.9s) | Q8 FEM with 2,625 nodes (PE = 4.420E-01, CPU=2.2s) | Q8 FEM with 504,452 nodes (PE = 4.419E-01, CPU=31.9s) |
| Number of Epochs NN-PU used to reach the FEM PE levels | 2 | 11 | 35 |
| CPU for NN-PU to reach the FEM PE levels (s) | 3.2 | 3.6 | 4.4 |

**Table 6: Comparison between the NN-PU solution obtained with known "parent" basis set and FEM solutions with different spatial discretization for the plate with two holes problems**

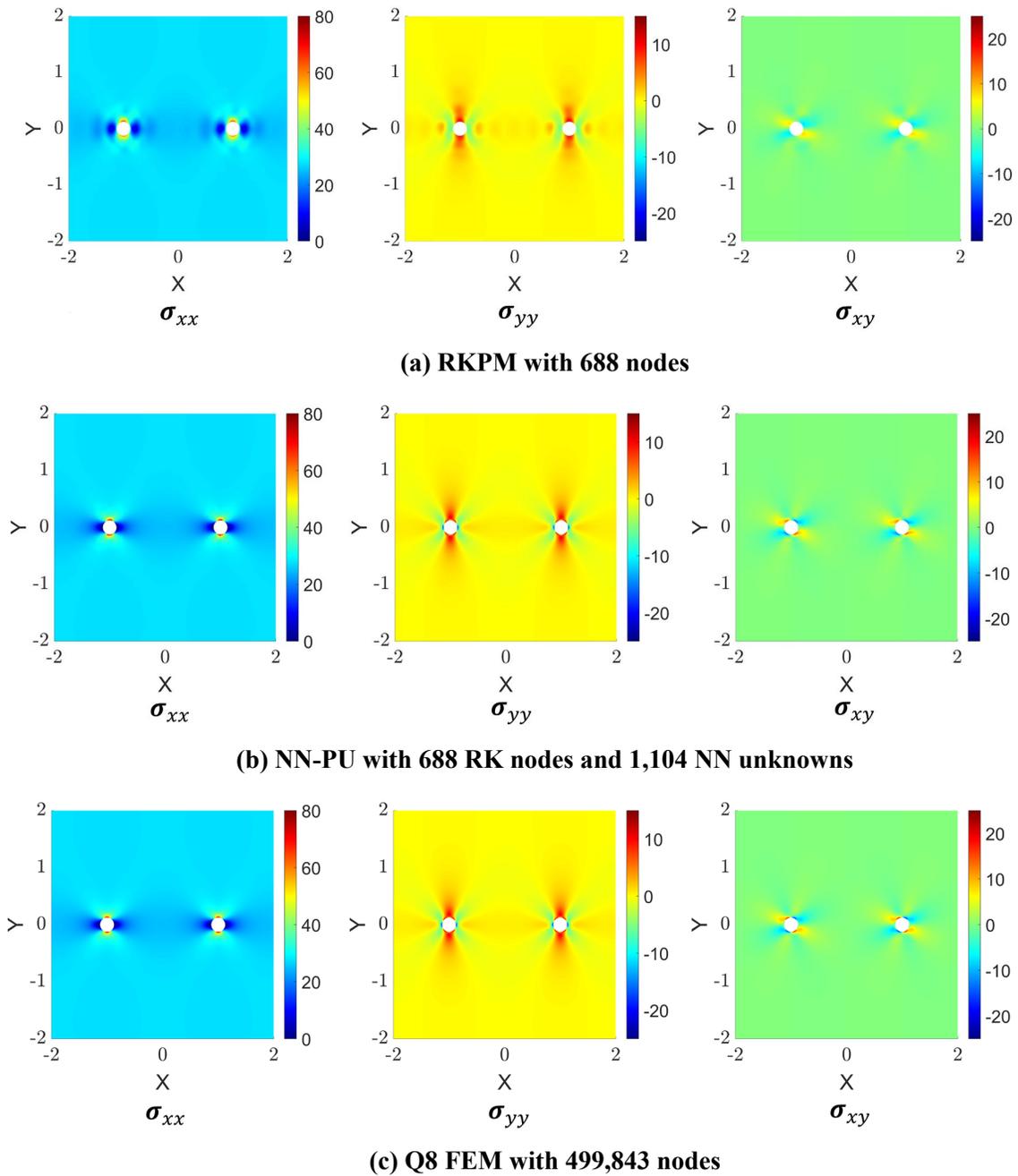

**(a) RKPM with 688 nodes**

**(b) NN-PU with 688 RK nodes and 1,104 NN unknowns**

**(c) Q8 FEM with 499,843 nodes**

Figure 36: Comparison of the stress solutions for plate with two small holes: (a) RKPM solutions with 688 nodes (b) NN-PU solutions with 688 RK nodes and 1,104 NN unknowns and (c) Q8 FEM solutions with 499,843 nodes

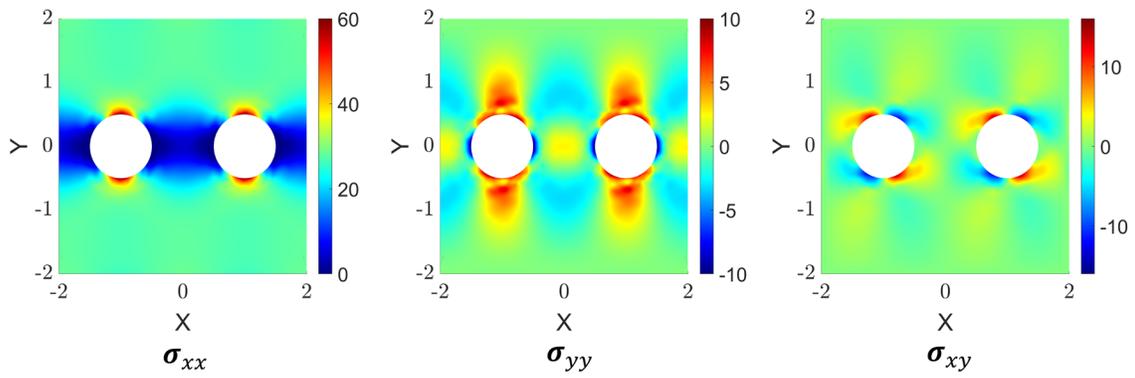

**(a) RKPM with 728 nodes**

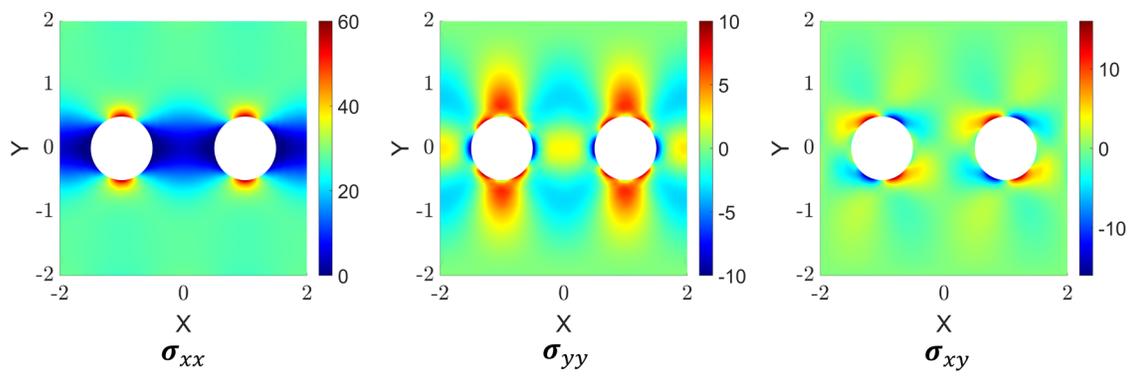

**(b) NN-PU with 728 RK nodes and 1,854 NN unknowns**

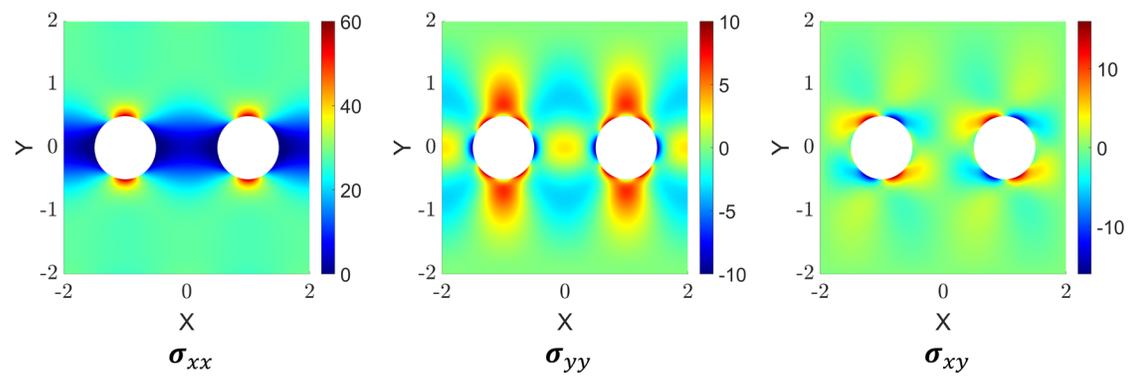

**(c) Q8 FEM with 504,452 nodes**

**Figure 37: Comparison of the stress solutions for plate with two big holes: (a) RKPM solutions with 728 nodes (b) NN-PU solutions with 728 RK nodes and 1,854 NN unknowns and (c) Q8 FEM solutions with 504,452 nodes**

### 5.3.2 Plate with multiple holes: transfer learned from the parent problems with single local features

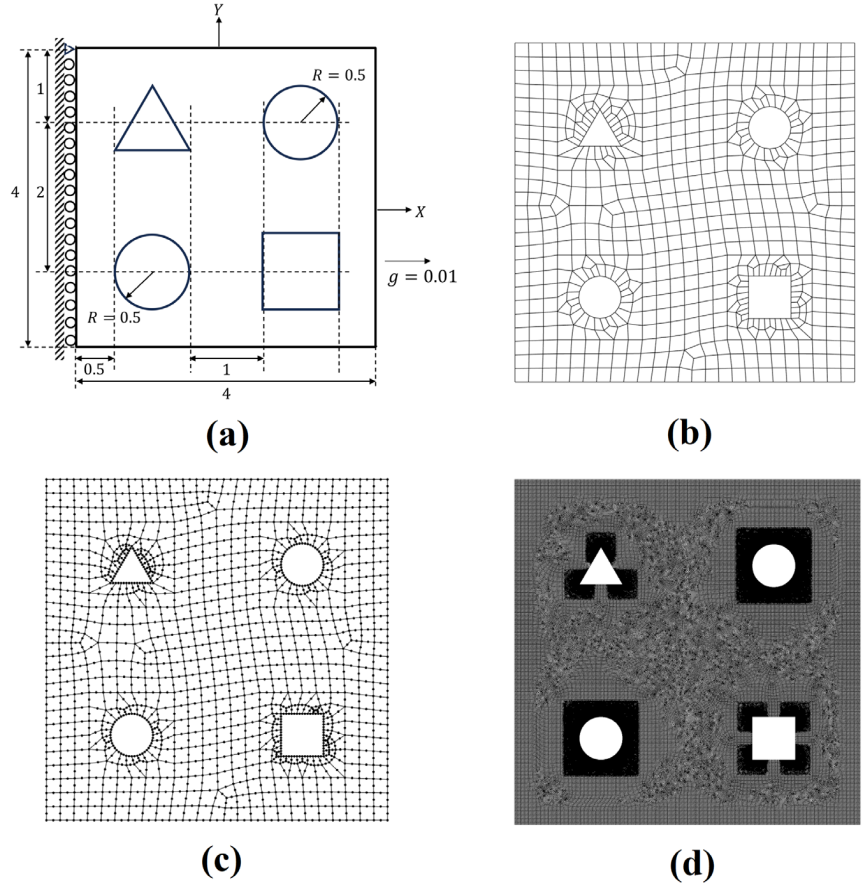

**Figure 38: (a) Problem setting for plate with multiple holes and 3 levels of FE discretization: (b) Q4 FEM with 768 nodes (c) Q8 FEM with 2,208 nodes (d) Q8 FEM with 647,174 nodes**

The last numerical example deals with a plate with multiple holes undergoing uniaxial tensile loading in the x-direction as depicted in Figure 38. In this case, 4 NN blocks are employed, with each block dedicated to capturing the localized characteristics around each hole. For these 4 NN blocks, 4 corresponding Parametric sub-blocks $\mathcal{N}_1^P$, $\mathcal{N}_2^P$, $\mathcal{N}_3^P$, and $\mathcal{N}_4^P$ are employed. All Parametric sub-blocks except the first one associated with the triangular hole utilize the same 2-layered architecture as that used in Sec. 5.2.1, and $\mathcal{N}_1^P$ adopts the same architecture as the one in Sec. 5.2.2. The NN-enrichment basis sets for each hole are taken from their corresponding "parent" problems in Sec. 5.1 with transfer learning. The NN Basis sub-block $\mathcal{N}_1^\zeta$, corresponding to the first NN block and associated with the triangular hole, are transfer-learned from the 90° "parent" L-shaped panel problem discussed in Sec. 5.1.2. Similarly, the NN Basis sub-blocks $\mathcal{N}_2^\zeta$ and $\mathcal{N}_3^\zeta$, corresponding to the two circular holes, are transfer-learned from the plate-with-a-hole "parent" problem in Sec. 5.1.1. The NN Basis sub-block $\mathcal{N}_4^\zeta$ corresponding to the square hole employs an eight-layered architecture with 20 neurons in the first 7 hidden layers and 5 NN-enrichment bases, which are learned from the 90° L-shaped panel "parent" problem discussed in Sec. 5.1.2. During the online calculation, all weights and biases in $\mathcal{N}_1^\zeta$, $\mathcal{N}_2^\zeta$, $\mathcal{N}_3^\zeta$, and $\mathcal{N}_4^\zeta$ are fully transferred from parent NN-enrichment basis sets learned from their corresponding "parent" problems. Consequently, the loss function minimization process only recomputes the optimal parametric coordinates and NN-enrichment coefficients within each NN block.

The NN-enriched nodesets for the 4 NN blocks are initiated as follows:

$$\mathcal{S}_1^\zeta = \{I | x_{1I} < 0, x_{2I} > 0, \rho_I^* > 0.1\rho_{max}^*\},$$
$$\mathcal{S}_2^\zeta = \{I | x_{1I} \geq 0, x_{2I} > 0, \rho_I^* > 0.1\rho_{max}^*\},$$
$$\mathcal{S}_3^\zeta = \{I | x_{1I} < 0, x_{2I} \leq 0, \rho_I^* > 0.1\rho_{max}^*\},$$
$$\mathcal{S}_4^\zeta = \{I | x_{1I} \geq 0, x_{2I} \leq 0, \rho_I^* > 0.1\rho_{max}^*\},$$
(40)

where $x_{1I}$ and $x_{2I}$ denote the x and y directional coordinates of a node $I$, and the subscripts in the NN-enriched nodesets are associated with the indices of the NN blocks. The assessment of nodal energy error densities takes place every 20 epochs, and the 4 enrichment node subsets are updated independently, as demonstrated in Figure 39, resulting a total of 2,990 trainable NN-enrichment coefficients.

The loss function minimization in the online calculation is illustrated in Figure 40. Table 7 reports a performance between NN-PU and FEM methods with varying degrees of domain discretization.

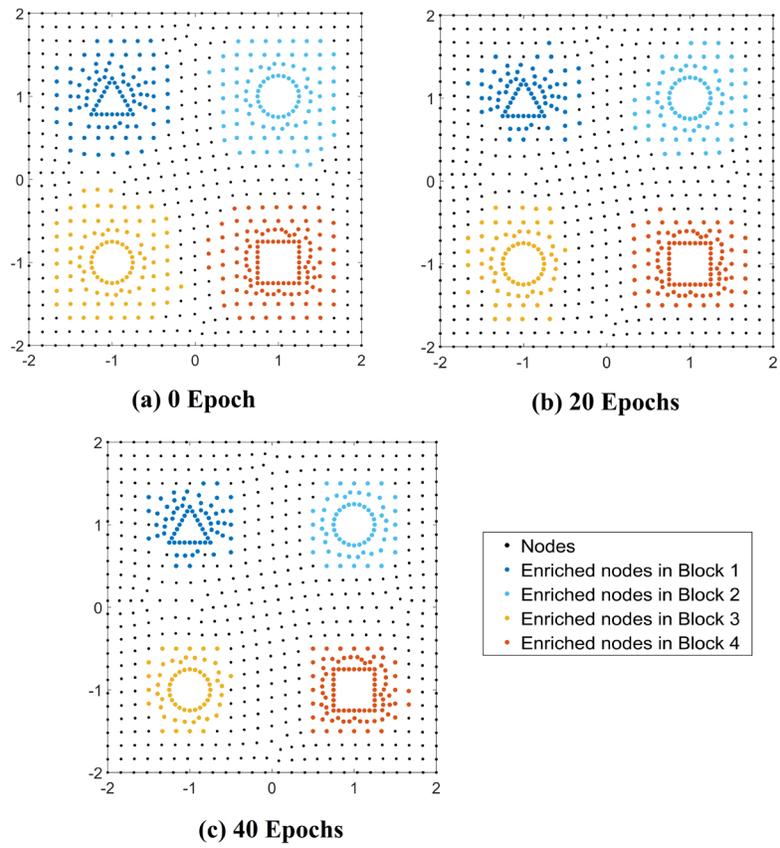

**(a) 0 Epoch**

**(b) 20 Epochs**

**(c) 40 Epochs**

**Figure 39: Adaptive selections of the NN-enriched nodes for 4 NN blocks: (a) total 442 nodes are enriched (b) total 351 nodes are enriched and (c) total 299 nodes are enriched** (numbers of enriched nodes in x- and y- directions are the same)

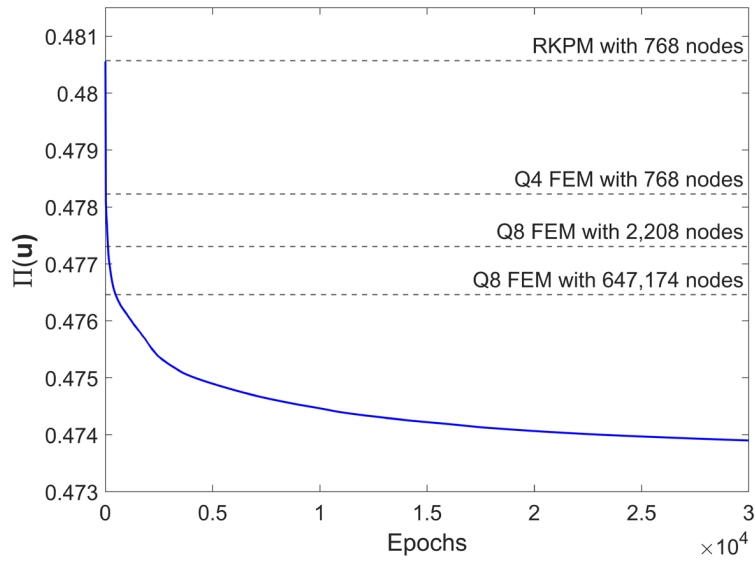

**Figure 40: Loss function minimization for NN-PU and the comparison with FEM potential energy**

| NN-PU with potential energy (PE) lower than FEM PE | Q4 FEM with 768 nodes (PE = 4.782E-01, CPU=1.9s) | Q8 FEM with 2,208 nodes (PE = 4.773E-01, CPU=2.2s) | Q8 FEM with 647,174 nodes (PE = 4.765E-01, CPU=45.2s) |
|---|---|---|---|
| Number of Epochs NN-PU used to reach the FEM PE levels | 17 | 120 | 480 |
| CPU (s) for NN-PU to reach the FEM PE levels | 3.9 | 8.1 | 22.4 |

**Table 7: Comparison between the NN-PU solution obtained with known NN enrichment basis sets and FEM solutions with different spatial discretization for the plate with four holes problem**

The stress results comparison between RKPM, NN-PU and the high-order highly refined FEM with over a million degrees of freedom is demonstrated in Figure 41. The NN-PU approximation with transfer-learned NN-enrichment basis sets is able to achieve solutions with similar accuracy to the much-refined FEM with 50.4% CPU reduction. The online calculation of NN-PU yields a converged potential energy of 4.7391E-01, resulting in a 0.5% potential energy reduction compared to that computed from Q8 FEM solutions with 647,174 nodes. Note that other than utilizing the NN enrichment basis sets pre-trained from "parent" problems in Sec. 5.1, $\mathcal{N}_1^\zeta$ can be transfer learned from the converged NN bases in the 60° L-shaped panel problem in Sec. 5.2.2. In this way, $\mathcal{N}_1^P$ can take the same 2-layered architecture as the other NN blocks, totaling of 3,098 NN unknowns and leading to a larger CPU reduction (61.5%) compared to the much-refined FEM. This study demonstrates that the NN-PU approach with transfer learning from parent problems with different features becomes more effective than the FEM approach when demanding higher solution accuracy.

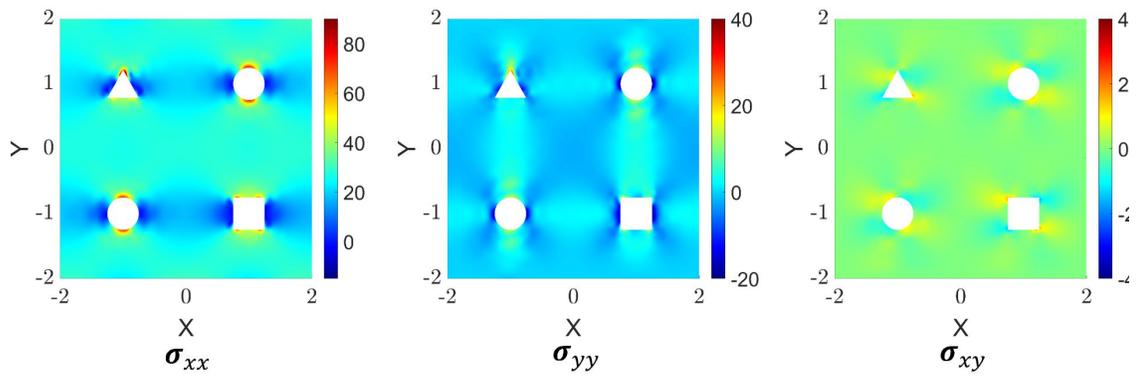

**(a) RKPM with 768 nodes**

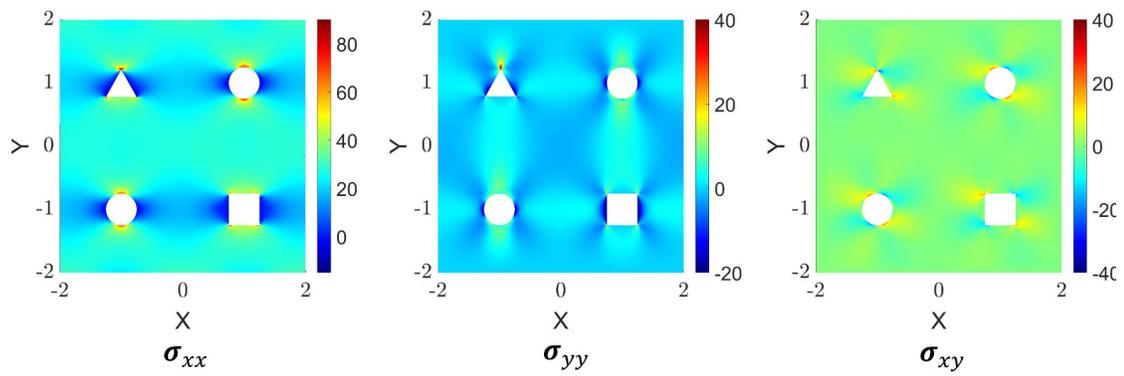

**(b) NN-PU with 768 RK nodes and 3,593 NN unknowns**

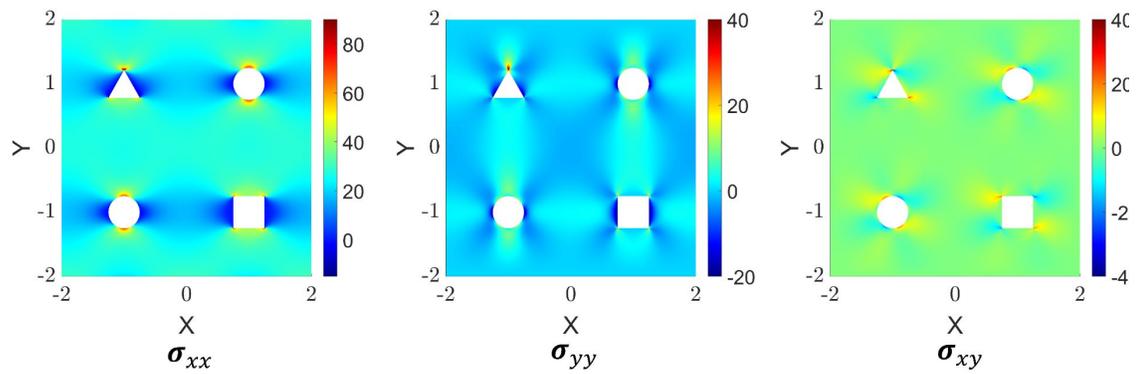

**(c) Q8 FEM with 647,174 nodes**

**Figure 41: Comparison of the stress solutions: (a) RKPM solutions with 768 nodes (b) NN-PU solutions with 768 RK nodes and 3,593 NN unknowns and (c) Q8 FEM solutions with 647,174 nodes**

# 6    Conclusion

A neural network-enriched Partition of Unity (NN-PU) method has been proposed for solving boundary value problems. In this approach, the domain geometry, boundary conditions and loading conditions are discretized by the background approximation that possesses partition of unity properties. The background approximation with a fixed coarse discretization is enriched adaptively by feature-encoded local NN basis functions trained at the offline stage via artificial neural networks with minimization of a potential energy-based loss function. The transferable block-level NN architecture used for the NN approximation keeps high sparsity of the network structure and offers flexibility for integrating the pre-trained feature-encoded NN enrichment basis sets with the background PU approximation, facilitating a more efficient online solution strategy. This block neural network has the flexibility to incorporate multiple feature-encoded NN bases targeting different localized features in the solution. The proposed NN-PU method, coupled with feature-encoded transfer learning, provides an effective adaptivity framework for solving PDEs, significantly simplifying the tedious conventional mesh-based h- and p-refinement procedures subjected to conforming constraints. Error analysis for the proposed NN-PU with reproducing kernel as the background PU approximation showing how the convergence is determined by both background PU approximation as well as the NN enrichment.

Employing a coarse background discretization, the method utilizes extrinsic NN-based enrichment functions to enhance the background approximation within the Partition of Unity framework. To improve the efficiency and adaptivity of the NN-PU approach, a block-level NN approximation is introduced, where each NN block is designed to target a specific underlying feature in the problem's physics. The architecture of the block-level NN approximation is designed with three sub-blocks: (1) the Parametric Sub-block to facilitate the utilization of pre-trained NN bases in a parametric coordinate to represent functions with local features, (2) the NN Basis Sub-block taking the parametric coordinates to generate NN enrichment basis functions and to embed the trained NN bases to the background approximation functions, (3) the Coefficient Sub-block to determine coefficients for linear combination of NN enrichment bases by optimizing an energy-based loss function and add the contribution of each NN enrichment basis to the total solution of the problem. During offline training, multiple NN enrichment basis sets are pre-trained for different "parent" problems to embed various underlying solution features into the NN basis functions. These feature-encoded NN basis sets are then utilized in online computation stages through transfer learning to capture complex solution patterns efficiently.

Parametric studies on the effects of different NN architectures in minimizing NN loss functions show that employing deeper and wider neural networks (up to 50 neurons in hidden layers) and increasing the number of NN bases (up to 5) significantly improve the speed of offline training of NN enrichment basis sets. An error indicator is introduced to identify the initial enrichment nodes, and the enrichment nodes are continuously updated during the loss function minimization iteration to facilitate the adaptive NN enrichment. The effectiveness of the NN-PU method in solving elasticity problems with local features has been substantiated through numerical examples. Comparison with FEM solutions obtained via Galerkin approximation has been made to examine the effectiveness of the proposed method. Numerical error analysis of a "parent" plate-with-a-hole problem shows good agreement on the analytically derived convergence rates. Assessments conducted during online solution calculations for problems featuring single or multiple local features indicate that the NN-PU method, incorporating feature-encoded transfer learning with pre-trained NN-enrichment basis sets, achieves comparable solution accuracy to highly refined FEM with considerably reduced unknowns and CPUs. This NN-PU framework can be extended to inelastic problems by a proper choice of energy function.

## Acknowledgements

The support of this work by the Sandia National Laboratories to UC San Diego under Contract Agreement 1655264 is greatly acknowledged. Sandia National Laboratories is a multi-mission laboratory managed and operated by National Technology and Engineering Solutions of Sandia, LLC, a wholly owned subsidiary of Honeywell International Inc., for the U.S. Department of Energy's National Nuclear Security Administration under contract DE-NA0003525. This paper describes objective technical results and analysis. Any subjective views or opinions that might be expressed in the paper do not necessarily represent the views of the U.S. Department of Energy or the United States Government. The support from National Science Foundation under Award Number CMMI-1826221 to University of California, San Diego is greatly appreciated.